\documentclass[a4paper,10pt]{article}
\usepackage[utf8]{inputenc}
\usepackage{amsmath}
\usepackage{amsfonts}
\usepackage{amssymb}
\usepackage{indentfirst}
\usepackage[english]{babel}
\usepackage{graphicx}
\usepackage{amscd}
\usepackage{tikz}
\usepackage{amssymb}
\usepackage{fancyhdr}
\newtheorem{defn}{Definition}[section]
\newtheorem{thm}[defn]{Theorem}
\newtheorem{cor}[defn]{Corollary}
\newtheorem{lem}[defn]{Lemma}
\newtheorem{rmk}[defn]{Remark}

\newtheorem{prop}[defn]{Proposition}

\newtheorem{nota}[defn]{Notations}
\newtheorem{pro}{Proposition}
\newtheorem{theo}{Theorem}

\newtheorem{coro}{Corollary}
\title{ \Large{\textbf{Mobile Product and Zariski decomposition}}\\
\vspace*{5mm}
\small{\textbf{MARIO PRINCIPATO}\\
\begin{center}
Universit\`a degli studi di Roma ``Tor Vergata'',
Via della Ricerca Scientifica 1, 0133 Rome, Italy\\
\emph{email:}\textit{principa@mat.uniroma2.it} \\
phone-number: +39 06 72594648
\end{center}}}
\begin{document}

\maketitle

\begin{abstract}
\flushleft We explain the relationship between $\alpha_{1}\cdots\alpha_{q}$ (standard cohomology product) and $(\alpha_{1}\cdots\alpha_{q})$ (mobile
intersection product) of pseudo-effective classes $\alpha_{1},\dots,\alpha_{q}$ on a compact K$\ddot{\text{a}}$hler manifold.  We also
show how to use this relationship for proving some holomorphic Morse inequalities. Then we prove a result concerning the direct image of
Lelong numbers under a modification in dimension 3,  deriving a continuity property for the Lelong numbers of the wedge of $(1,1)-$currents.
\end{abstract}

\begin{center}\section*{Introduction}\end{center}
\flushleft
 Since the  work of O. Zariski \cite{Zar} the study of the ring
\begin{equation*}
 \displaystyle R(X,D):=\bigoplus_{k\geq0}H^{0}(X,\mathcal{O}(kD)),
\end{equation*}
where $D$ is an effective divisor on a projective surface X, became very important.\\
In particular the main result  is that any $\mathbb{Q}$-divisor $D$ on a projective surface $X$ can be  decomposed into
a sum $D=P+N$ where $P$ is a nef $\mathbb{Q}$-divisor, $N=\sum a_{j}D_{j}$ is an effective $\mathbb{Q}$-divisor
such that  $(D_{i}\cdot D_{j})$ is negative definite, and $P$ is orthogonal to $N$ with respect to the intersection form.
Then Zariski showed that  $H^{0}(kP)\hookrightarrow H^{0}(kD)$ is  an isomorphism in that case.\\
Using  the  metrics with minimal singularities on a pseudo-effective line bundle $L$ introduced by J.-P. Demailly,
S. Boucksom \cite{Bou1} defines the divisorial Zariski decomposition of a pseudo-effective class $\alpha=(\alpha)+N(\alpha)$ where $N(\alpha)$ is an 
effective $\mathbb{R}$-divisor which is ``exceptional`` in some sense.\\
In this work we are going to study the relationship between the product $\alpha^{i}$ and $(\alpha^{i})$  with respect the
\emph{non-nef locus} of $\alpha$ $E_{nn}(\alpha)$ as defined in (\cite{Bou1}),  where $(\alpha^{i})$
is the \textit{mobile intersection product} introduced in (\cite{BDPP}).
 Let $\mathcal{E}$ be the cone of pseudo-effective classes, we can prove
\begin{pro}
 Let $(X,\omega)$ be a   compact complex K$\ddot{\text{a}}$hler manifold of complex
dimension $n$ and let $\alpha_{1},\dots,\alpha_{q}\in\mathcal{E}$ such that $cod(E_{nn}(\alpha_{1})\cup\dots\cup E_{nn}(\alpha_{q}))=m$.\\ 
Then if $s<m$ and $j_1<j_2<\cdots<j_s\in\{1,\cdots,q\}$ we have
\begin{equation}
 \alpha_{j_{1}}\cdots\alpha_{j_{s}}=(\alpha_{j_{1}}\cdots\alpha_{j_{s}})
\end{equation}

 In particular when $\alpha=\alpha_{1}=\dots=\alpha_{q}$ and $cod(E_{nn}(\alpha))=q$ then 
\begin{equation}
 \alpha^{j}=\langle\alpha^{j}\rangle=(\alpha^{j})\quad\forall j=1,\dots,q-1.
\end{equation}
\end{pro}
Furthermore we can study the extremal case of Proposition above, and introducing for $0\leq q\leq n$ the cone
\begin{equation*}
 \mathcal{M}^{q}:=\{\alpha\in\mathcal{E}:~\dim(E_{nn}(\alpha))\leq n-q\},
\end{equation*} 
then we have the following
\begin{pro}
 $\mathcal{M}^{q}\subseteq\mathcal{E}$ is a closed convex cone.
\end{pro}
Moreover one has
\begin{pro}
 $E_{nn}(\alpha)$ doesn't have irreducible components of dimension zero.
\end{pro}
Thus we can prove the following
\begin{theo}
 Let $X$ be a compact K$\ddot{\text{a}}$hler manifold and be $\alpha_{1},\dots,\alpha_{q}\in{\mathcal{M}}^{q}$.
Let $\{Y_{h}\}_{h\in H}$ the family of codimension $q$ components of $\displaystyle\bigcap_{j=1}^{q}E_{nn}(\alpha_{j})$.
Then we have the following decomposition:\\

\item i) $\displaystyle\alpha_{1}\cdots\alpha_{q}=\langle\alpha_{1}\cdots\alpha_{q}\rangle+\left\{\sum_{h\in H}\nu(\alpha_{1}\cdots\alpha_{q},Y_{h})[Y_{h}]\right\},~\nu(\alpha_{1}\cdots\alpha_{q},Y_{h})>0$
\\\hspace*{2mm} where $\nu(\alpha_{1}\cdots\alpha_{q},Y_{h})$ are \hspace*{2mm} the multiplicities of $\alpha_{1}\cdots\alpha_{q}$ (see Definition \ref{defini}),
and the set $H$ is at most countable;
\\
\item ii) if $\alpha_{j}\in \text{int}(\mathcal{M}^{q})$ for $1\leq j\leq q$, then the set $H$ is finite and\\ \hspace*{2mm}
$\nu(\alpha_{1}\cdots\alpha_{q},Y_{h})=\nu(T_{min,\alpha_{1}}\wedge\cdots\wedge T_{min,\alpha_{q}},Y_{h})$ for all $h\in H$,
where $T_{min,\alpha_{j}}$ is a positive current with minimal singularities in $\alpha_{j}$.
\end{theo}
As an immediate consequence of \textbf{Theorem 1} we have
\begin{coro}
 If $\alpha\in\emph{int}(\mathcal{M}^{q})$ then $\nu(T_{min}^{q},Y)=0$ for every irreducible analytic set of codimension $q$ not completely
contained in $E_{nn}(\alpha)$.
\end{coro}

We can give a partial converse of this Theorem (\ref{thm3} (ii))
\begin{pro}
 Let $(X,\omega)$ be a compact K$\ddot{\text{a}}$hler manifold of complex dimension $n$, and let $\alpha\in\mathcal{E}$. Then
$\alpha\in\mathcal{M}^{q+1}$ if and only if $\alpha^{j}=(\alpha^{j})$ for all $j=1,\dots,q$.\\
In particular $\alpha$ is nef if and only if $\alpha^{j}=(\alpha^{j})$ for all $j=1,\dots,n-1$.
\end{pro}

Then we focused our attention to the Holomorphic Morse Inequalities, a theory initiated by J.-P Demailly in the '80's.\\
Let's suppose we have an hermitian line bundle $L$ over a compact K$\ddot{\text{a}}$hler manifold and we want to study the asymptotic behavior (for $k\to+\infty$)
of the partial alternating sum of $h^{q}(X,kL)$; the complete sum is simply the Euler characteristic. In general the behavior is controlled by an
estimate involving the integral of the top wedge power of the Chern curvature of $L$, extended over suitable subsets of $X$.
One difficulty in the application of these inequalities is that the curvature integral is in general quite uneasy to compute, it is neither a topological
nor an algebraic invariant. However a special case of the Morse inequalities can be reformulated in a more algebraic setting in which only algebraic invariants
are involved, see e.g. \cite{Trap1}, \cite{Siu2}.\\
Now by using  the decomposition
\begin{equation*}
 \alpha_{1}\cdots\alpha_{q}-\langle\alpha_{1}\cdots\alpha_{q}\rangle=\left\{\sum_{h=1}^{t}\nu(\alpha_{1}\cdots\alpha_{q},Y_{h})[Y_{h}]\right\}
\end{equation*}
we can prove
\begin{theo}
Let $X$ be a compact projective manifold of complex dimension $n$. Let $L$ and $F$ be two line bundles over $X$ with $L$ nef and
$F\in\mathcal{M}^{s}$ with $\dim (E_{nn}(F))=n-s$. Let $\{Y_{t}\}_{t\in T}$ be the irreducible components (possibly infinite) of codimension s of 
$E_{nn}(\alpha)$, and
let $\nu_{t}$ and $\nu_{t}^{'}$  be the multiplicities of $F^{s}$ and the multiplicity of $F$ along
$Y_{t}$ respectively. Then we have the following \emph{strong} Morse inequalities
\begin{multline*}
 \displaystyle \sum_{j=0}^{s}(-1)^{s-j}h^{j}(X,k(L-F))\leq
 \frac{k^{n}}{n!}(\sum_{j=0}^{s}\binom{n}{j}(-1)^{s-j}L^{n-j}(F^{j})+\\
 \sum_{t=1}^{+\infty}\binom{n}{s}(L+\nu_{t}^{'}\{u\})^{n-s}\nu_{t}[Y_{t}])+o(k^n).
\end{multline*}
\end{theo}
And as a consequence we have
\begin{coro}
 Let $X$ be a compact projective manifold of complex dimension $n$. Let $L$ and $F$ be two line bundles over $X$ with $L$ nef and
$F\in\mathcal{M}^{s}$ with $\dim (E_{nn}(F))=n-s$. Then we have the following Morse inequalities
\begin{multline*}
 \displaystyle \sum_{j=0}^{s}(-1)^{s-j}h^{j}(X,k(L-F))\leq\\
 \frac{k^{n}}{n!}\sum_{j=0}^{s}\binom{n}{j}(-1)^{s-j}L^{n-j}(F^{j})+\\
 \binom{n}{s}(L+b\{u\})^{n-s}(F^{s}-\langle F^{s}\rangle)+o(k^n),
\end{multline*}
where $\displaystyle b:=\max_{t=1,\dots,N}\nu_{t}$.\end{coro}
In the last Section we restrict to the dimension 3 case and we prove the following result about  the direct image of Lelong numbers under a modification,
\begin{theo}
Let $X$ be a  complex compact manifold with \\ $\dim_{\mathbb{C}}(X)=3$.
Let $\tilde{\mu}:\tilde{X}\rightarrow X$ be a modification of $X$ and $\Omega$ is a smooth, positive form on $\tilde{X}$ of bidimension 1. Then
$\nu(\tilde{\mu}_{*}(\tilde{\Omega}),Y)=0$ $\forall Y$  irreducible curve on $X$, where $\nu$ is the generic Lelong number.
\end{theo}
Thus we can prove the following continuity property for Lelong numbers
\begin{pro}
 Let $\alpha,\beta\in\emph{int}(\mathcal{M}^{2})$. Then
\begin{equation*}
 \displaystyle\lim_{k\to+\infty}\nu(T_{k,\alpha}\wedge T_{k,\beta},x)=\nu(T_{min,\alpha}\wedge T_{min,\beta},x)~\forall~x\in X,
\end{equation*}
where $\{T_{k,\alpha}\}$ and $\{T_{k,\beta}\}$ are sequences of currents with analytic singularities obtained in \cite{Dem1} which converge weakly to $T_{min,\alpha}$
and $T_{min,\beta}$ respectively. 
\end{pro}

\
\newline
\textbf{Acknowledgements}\\ I would like thank my advisor Stefano Trapani for introducing me to the
topic of this article, for his great knowledge and inspiration, and for his continuous support.
\newpage
\section{Decomposition(s) in $H_{\geq0}^{k,k}(X)$}
We  need two basic types of regularizations (inside a fixed cohomology class) for $(1,1)-$current, both due to J.-P. Demailly \cite{Dem1} 
\begin{thm}\label{thm2}
 Let $T=\alpha+i\partial\bar{\partial}\varphi$ be a closed almost positive $(1,1)-$current on a compact complex manifold $X$ and fix an Hermitian form $\omega.$ Suppose that
$T\geq\gamma$ for some smooth real $(1,1)-$form $\gamma$ on $X$. Then:\\
i) There exists a sequence of smooth forms $\theta_{k}$ in $\{T\}$ (cohomology class of T) which converges weakly to $T$ and such that 
$\theta_{k}\geq\gamma-C\lambda_{k}\omega$ where $C$ is a constant depending on the curvature of $(T_{X},\omega)$ only, and $\lambda_{k}$ is a decreasing
sequence of continuous functions such that $\lambda_{k}(x)\to\nu(T,x)$ for every $x\in X.$\\
ii) There exists a sequence $T_{k}=\alpha+i\partial\bar{\partial}\varphi_{k}$ of closed currents such that:
\begin{itemize}
 \item $\varphi_{k}$ (and thus $T_{k}$) is smooth on the complement $X\setminus Z_{k}$ of an analytic set $Z_{k}$ such that there is an increasing sequence
\begin{equation*}
 Z_{0}\subset Z_{1}\subset\dots\subset X.
\end{equation*}
\item There is a uniform estimate $T_{k}\geq\gamma-\varepsilon_{k}\omega$ with $\lim\downarrow\varepsilon_{k}=0$ as $k$ tends to $+\infty$.\\
\item The sequence $(\varphi_{k})$ is non increasing, and we have $\lim\downarrow\varphi_{k}=\varphi$. As a consequence, $T_{k}\rightharpoonup T$.\\
\item Near $Z_{k}$, the potential $\varphi_{k}$ has logarithmic poles,namely, for every $x_{0}\in Z_{k}$, there is a neighborhood $U$
of $x_{0}$ such that $\varphi_{k}(z)=\lambda_{k}\log\sum_{l}|g_{k,l}|^{2}+O(1)$ for suitable holomorphic functions $(g_{k,l})$ on $U$ and $\lambda_{k}>0$.
Moreover, there is a (global) proper modification $\mu_{k}:X_{k}\to X$ of $X$, obtained as a sequence of blow-ups with smooth
centers, such that $\varphi_{k}\circ\mu_{k}$ can be written locally on $X_{k}$ as
\begin{equation*}
 \varphi_{k}\circ\mu_{k}(w)=\lambda_{k}\left(\sum n_{l}|\tilde{g}_{l}|^{2}+f(w)\right)
\end{equation*}
where $(\tilde{g}_{l}=0)$ are local generators of suitable (global) divisor $D_{l}$ on $X_{k}$
such that $\sum D_{l}$ has normal crossings, $n_l$ are positive integers and the $f$'s are smooth functions on $X_{k}.$
\end{itemize}
 \end{thm}
It is important to observe that such $\mu_{k}$ are obtained blowing up along the  multiplier ideal sheaves $\mathcal{I}(k\varphi)$ defined as
\begin{equation*}
 \mathcal{I}(k\varphi)_{x}:=\left\{f\in\mathcal{O}_{\Omega,x}:~|f|^{2}e^{-2k\varphi}\in\text{L}_{loc}^{1}\right\}
\end{equation*}
where $\Omega\subseteq X$ is open.
\newpage

\flushleft Following  \cite{Bou1}, when $\alpha\in\mathcal{E}$ one can introduce a measure of nefness of $\alpha$.\\
Let $\psi_{1},\psi_{2}$ be almost plurisubharmonic functions on $X$ we say that $\psi_1$ is less
singular than $\psi_2$ (and write $\psi_{1}\preceq\psi_{2}$) if $\psi_{2}\leq\psi_{1}+C$ for some constant $C$. When $S_{1}$ and $S_{2}$ are closed almost positive (1,1)-currents on $X$, 
we can  compare their singularities  by comparing those of their local potentials $\psi_{1},\psi_{2}$.\\
For each $\varepsilon>0$ let $T_{min,\varepsilon}$ be a current with minimal singularities in
 $\alpha[-\varepsilon\omega]$ which is the set of closed almost positive (1,1)-currents $T$ lying in $\alpha$ with $T\geq-\varepsilon\omega.$
\begin{defn}
 The minimal multiplicity at $x\in X$ of the pseudo-effective class $\alpha$ is defined as
\begin{equation*}
 \displaystyle\nu(\alpha,x)=\sup_{\varepsilon>0}~\nu(T_{min,\varepsilon},x).
\end{equation*}
\end{defn}
When $D$ is a prime divisor, we define the generic minimal multiplicity of $\alpha$ along $D$ as
\begin{equation*}
 \nu(\alpha,D)=\inf\{\nu(\alpha,x),x\in D\}.
\end{equation*}
We then have $\displaystyle\nu(\alpha,D)=\sup_{\varepsilon>0}\nu(T_{min,\varepsilon},D), $ and $\nu(\alpha,D)=\nu(\alpha,x)$ for the very
general $x\in D.$\\
One can give the following
\begin{defn}
 The non-nef locus of a pseudo-effective class $\alpha$ is defined by
\begin{equation*}
 E_{nn}(\alpha):=\{x\in X,~\nu(\alpha,x)>0\}.
\end{equation*}
\end{defn}
\begin{rmk}\label{rmk9}
 Let $\alpha$ be a pseudo-effective cohomology class and be $\varepsilon$ small positive rational number.  Let's fix a smooth hermitian
form $\omega$ on $X$ so that $T_{min,\varepsilon}+\varepsilon\omega$ is a positive current representing $\alpha+\varepsilon\omega$.
 Then 
$T_{min,\varepsilon}+\varepsilon\omega\succeq T_{min,\alpha+\varepsilon\omega}$ where $T_{min,\alpha+\varepsilon\omega}$ is a positive current with minimal singularities
in $\alpha+\varepsilon\omega$. On the other hand $T_{min,\alpha+\varepsilon\omega}-\varepsilon\omega$ is a current representing $\alpha$ such that
$T_{min,\alpha+\varepsilon\omega}-\varepsilon\omega\geq-\varepsilon\omega$ then 
$T_{min,\alpha+\varepsilon\omega}-\varepsilon\omega\succeq T_{min,\varepsilon}$. Hence we infer that $T_{min,\varepsilon}+\varepsilon\omega$
is a positive current with minimal singularities representing $\alpha+\varepsilon\omega$. 
\end{rmk}

Now let $X$ be a compact K$\ddot{\text{a}}$hler manifold and let $Y\subseteq X$ be an analytic set of dimension $p$ and let
$\alpha_{1},\dots,\alpha_{p}\in H^{1,1}(X,\mathbb{R})$ be pseudo-effective classes. Then one can define a ``positive'' number of intersection
$(\alpha_{1}\cdot\cdots\cdot\alpha_{p}\cdot Y)_{>0}$ which intuitively should be equal to  the usual intersection number $(\beta_{1}\cdot\cdots\cdot\beta_{p}\cdot Y)$,
where $\beta_i$ is the nef part of $\alpha_i$ in its Zariski decomposition. But in general the Zariski decomposition does not exist, however
using the currents with analytic singularities, one can solve the problem. More generally  one can consider a closed positive current $\Theta$
instead of $Y$, then one can give the following
\begin{defn}[\cite{Bou1}]\label{def3}
 The mobile intersection number $(\alpha_{1}\cdot\cdots\cdot\alpha_{p}\cdot \Theta)_{>0}$ of $\alpha_i$'s and $\Theta$ is defined as
\begin{equation*}
 \displaystyle\inf_{\varepsilon>0}\left(\sup\int_{X\setminus F}(T_{1}+\varepsilon\omega)\wedge\dots\wedge(T_{p}+\varepsilon\omega)\wedge\Theta\right),
\end{equation*}
where $T_{i}\in\alpha_{i}[-\varepsilon\omega]$ are currents with analytic singularities and $F$ is the union of its unbounded-loci.
\end{defn}

Let's observe that the integrals in the definition above are all convergent and the supremum is finite
because the integrals can be bounded in terms  only of cohomology classes of currents, moreover the definition does not depend on the choice of 
the K$\ddot{\text{a}}$hler form $\omega$. Finally the supremum increases with $\varepsilon$
so that the limit exist and it is equal to the infimum for $\varepsilon>0$. 
Then by duality we have
\begin{thm}[\cite{BDPP}]\label{thm3}
 Let $X$ be a compact K$\ddot{a}$hler manifold. We denote here by $H^{k,k}_{\geq0}(X)$ the cone of cohomology classes of type $(k,k)$ which
have non-negative intersection with all closed semi-positive smooth forms of bidegree $(n-k,n-k)$.
\item (i)For each $k=1,\dots,n,$ there exists a canonical  $\mathbf{mobile~ intersection~ product}$
\begin{equation*}
 \mathcal{E}\times\cdots\times\mathcal{E}\to H^{k,k}_{\geq0}(X),~(\alpha_{1},\dots,\alpha_{k})\to\langle\alpha_{1}\cdot\alpha_{2}\cdots\alpha_{k}\rangle
\end{equation*}

\item (ii) The product is increasing, homogeneous of degree 1 and superadditive in each argument, i.e.
\begin{equation*}
 \langle\alpha_{1}\cdots(\alpha_{j}^{'}+\alpha_{j}^{''})\cdots\alpha_{k}\rangle\geq\langle\alpha_{1}\cdots\alpha_{j}^{'}\cdots\alpha_{k}\rangle+
\langle\alpha_{1}\cdots\alpha_{j}^{''}\cdots\alpha_{k}\rangle.
\end{equation*}
It coincides with the ordinary intersection product when the $\alpha_{j}$ are nef classes.
\end{thm}

\begin{rmk}\label{rmk7}
 Let us note that by construction $(\alpha_{1}\cdots\alpha_{p})\cdot u=(\alpha_{1}\cdots\alpha_{p}\cdot u)_{>0}$ for all closed semi-positive
forms $u.$
\end{rmk}

Now using the definition of \emph{non-pluripolar product} $\langle u_{1}\cdots u_{p}\rangle$ (\cite{BPGZ})  for $u_{1},\dots,u_{p}$  psh-function on an open set of $X$ (here $X$ is a compact complex manifold),
one can define a cohomology class $\langle\alpha_{1}\cdots\alpha_{p}\rangle\in H^{p,p}(X,\mathbb{R})$ as follows
\begin{defn}
 Let $\alpha_{1},\dots,\alpha_{p}\in H^{1,1}(X,\mathbb{R})$ be big cohomology classes and let $T_{min,i}\in\alpha_{i}$ be a positive current
with minimal singularities. Then the cohomology class of the non-pluripolar product $\langle T_{min,1}\wedge\dots\wedge T_{min,p}\rangle$ is indipendent of the 
choice of $T_{min,i}\in\alpha_{i}$ with minimal singularities. It will be denoted by
\begin{equation*}
 \langle\alpha_{1}\cdots\alpha_{p}\rangle\in H^{p,p}(X,\mathbb{R})
\end{equation*}
and called the $\mathbf{non-pluripolar~ product}$ of the $\alpha_{j}$. If $\alpha_{1},\dots,\alpha_{p}\in\mathcal{E}$ one sets
\begin{equation*}
 \displaystyle\langle\alpha_{1}\cdots\alpha_{p}\rangle=\lim_{\varepsilon\to0}\langle(\alpha_{1}+\varepsilon\beta)\cdots(\alpha_{p}+\varepsilon\beta)\rangle
\end{equation*}
where $\beta$ is an arbitrary K$\ddot{a}$hler class, using the continuity of the non-pluripolar product.
\end{defn}

Now on a compact K$\ddot{\text{a}}$hler manifold there exist two possible products (different a priori) : the mobile intersection product and the positive product of  pseudo-effective cohomology classes;
actually, as it is remarked in \cite{BPGZ} with no proof, these products are equal.
For the reader convenience we give a proof of the following well-known 
\begin{lem}\label{lem3}
 Let $T$ be any closed positive $(p,p)$-current on $X$ compact complex K$\ddot{\text{a}}$hler manifold. Then the Lelong number $\nu(T,x)$ of
$T$ can be bounded by a constant depending only on the $\partial\bar{\partial}-$cohomology class of $T$.
\end{lem}
\textit{Proof.} Let $\omega$ be a K$\ddot{\text{a}}$hler form on $X$, one has by definition that $\nu(T,x)$ is (up to a constant depending on $\omega$
near $x$) the limit for $r\to0^{+}$ of
\begin{equation*}
 \nu(T,x,r):=\frac{(n-p)!}{(\pi r^{2})^{n-p}}\int_{B(x,r)}T\wedge\omega^{n-p},
\end{equation*}
known to be an increasing function of $r$. Thus if we choose $r_{0}$ small enough to ensure that each ball $B(x,r_{0})$ is contained in a coordinate
chart, we get $\nu(T,x)\leq\nu(T,x,r_{0})\leq C\int_{X}T\wedge\omega^{n-p}$, a quantity depending only on the
cohomology class of $T$. \hfill$\square$ 

\begin{prop}\label{mobplu}
 Let $X$ be a compact K$\ddot{a}$hler manifold and let $\alpha_{1},\dots,\alpha_{p}\in H^{1,1}(X,\mathbb{R})$ be pseudo-effective classes. Then
\begin{equation}
 (\alpha_{1}\cdots\alpha_{p})=\langle\alpha_{1}\cdots\alpha_{p}\rangle.
\end{equation}
\end{prop}
\textit{Proof.} To keep notations simple we assume $\alpha=\alpha_{1}=\dots=\alpha_{p}$. Let's suppose first that $\alpha$ is big.
By remark (\ref{rmk7}) this is equivalent to showing that:
\begin{equation*}
 (\alpha^{p}\cdot u)_{>0}=\int_{X}\langle T_{min}^{p}\rangle\wedge u
\end{equation*}
for all closed semi-positive smooth forms $u$ of bidegree $n-p$, where $T_{min}$ is a positive current with minimal singularities in $\alpha.$
Let $T_{k}\in\alpha[-\varepsilon_k\omega]$ be a sequence of currents with analytic singularities such that $T_{k}\rightharpoonup T_{min}$.
Set $S_k:=(1-\varepsilon_k)T_k+\varepsilon_k S$ where $S$ is a K$\ddot{\text{a}}$hler current with analytic singularities such that 
$E_{+}(S)=E_{nK}(\alpha)$,
where
\begin{equation*}
 \displaystyle E_{nK}(\alpha):=\bigcap_{T\geq\varepsilon\omega,T\in\alpha}E_{+}(T)
\end{equation*}
  is the non-K$\ddot{\text{a}}$hler locus of $\alpha$ as in \cite{Bou1} and $E_{+}(S)=\bigcup_{c>0}\{x\in X~:~\nu(S,x)\geq c\}$. Then $\langle S_{k}^{p}\rangle\rightharpoonup\langle T_{min}^{p}\rangle$
on $X$ from which we infer 
\begin{equation*}
 \int_{X\setminus F_k}S_{k}^{p}\wedge u\leq(\alpha^{p}\cdot u)_{>0}
\end{equation*}
and letting $k\to+\infty$
\begin{equation*}
 \int_{X}\langle T_{min}^{p}\rangle\wedge u\leq(\alpha^{p}\cdot u)_{>0}.
\end{equation*}
But by the proof of theorem (\ref{thm3}) we can find a sequence of K$\ddot{\text{a}}$hler currents $R_k\in\alpha$ with analytic singularities such that
\begin{equation*}
 \int_{X\setminus G_k}R_{k}^{p}\wedge u\xrightarrow{k\to+\infty}(\alpha^{p}\cdot u)_{>0}
\end{equation*}
 where $G_k$ is the unbounded-locus of $R_k$. Then we find that
\begin{equation*}
 \int_{X\setminus G_k}R_{k}^{p}\wedge u\leq\int_{X}\langle T_{min}^{p}\rangle\wedge u
\end{equation*}
 and letting $k\to+\infty$ we find the other inequality.\\
If $\alpha_{1},\dots,\alpha_{p}$ are merely pseudo-effective since $(\alpha_{1}\cdots\alpha_{p})$ and $\langle\alpha_{1}\cdots\alpha_{p}\rangle$ depend
continuously on the $p$-tuple $\alpha_{1},\dots,\alpha_{p}$ of big classes the statement holds.\hfill$\square$
\\
\flushleft Now let $\alpha\in H^{1,1}(X,\mathbb{R})$ be a pseudo-effective class, we want to study the relation between $\alpha^{p}$ (the standard product in cohomology)
and $\langle\alpha^{p}\rangle=(\alpha^{p})$ when the codimension of $E_{nn}(\alpha)>p$. In general we have
\begin{prop}\label{pro2}
 Let $\alpha_{1},\dots,\alpha_{q}$ be pseudo-effective classes on compact complex K$\ddot{\text{a}}$hler manifold $(X,\omega)$ of complex
dimension $n$, such that $cod(E_{nn}(\alpha_{1})\cup\dots\cup E_{nn}(\alpha_{{q}}))=m$ 
for all $s_{1}<\dots<s_{t}\in\{1,\dots,q\}$. Then
\begin{equation}
 \alpha_{s_{1}}\cdots\alpha_{s_{t}}=(\alpha_{s_{1}}\cdots\alpha_{s_{t}})
\end{equation}
for all $t<m$.\\ In particular when $\alpha=\alpha_{1}=\dots=\alpha_{q}$ and $cod(E_{nn}(\alpha))=q$ then 
\begin{equation}
 \alpha^{j}=\langle\alpha^{j}\rangle=(\alpha^{j})\quad\forall j=1,\dots,q-1.
\end{equation}
\end{prop}
\textit{Proof.} Let $\{T_{k,s_{i}}\}\in\alpha_{s_{i}}$ be a sequence of currents with analytic singularities as in Theorem (\ref{thm2})
such that $T_{k,s_{i}}\rightharpoonup T_{min,s_{i}}$ for $i=1,\dots,t$, where $T_{min,s_{i}}\in\alpha_{s_{i}}$ is a current with minimal singularities with local potentials
$\varphi_{min,s_{i}}$. 
Let's define for all $i=1,\dots,t$
\begin{equation*}
 S_{k,s_{i}}:=T_{k,s_{i}}+\varepsilon_{k,s_{i}}\omega\geq0
\end{equation*}
and let $\mu_{k}:X_{k}\to X$ a common log-resolution for $S_{k,s_{i}}$ obtained  by blowing-up along $V=V(\mathcal{I}(k\varphi_{min,s_{1}}))\cup\dots\cup V(\mathcal{I}(k\varphi_{min,s_{t}}))$,
such that for all $i=1,\dots,t$
\begin{equation*}
 \mu_{k}^{*}S_{k,s_{i}}=[E_{k,s_{i}}]+\beta_{k,s_{i}}
\end{equation*}
where $\beta_{k,s_{i}}\geq0$ and smooth. Then
\begin{equation*}
 S_{k,s_{1}}\wedge\dots\wedge S_{k,s_{t}}-(\mu_{k})_{*}(\beta_{k,s_{1}}\wedge\dots\wedge\beta_{k,s_{t}})=0~\text{on}~X\setminus V,
\end{equation*}
and by using the support theorem for currents we get
\begin{equation*}
\displaystyle S_{k,s_{1}}\wedge\dots\wedge S_{k,s_{t}}-(\mu_{k})_{*}(\beta_{k,s_{1}}\wedge\dots\wedge\beta_{k,s_{t}})=\sum_{h=1}^{+\infty}\lambda_{h,k}[A_{h,k}]
\end{equation*}
where $[A_{h,k}]$ are the irreducible components of codimension $t$ in $V$, but $V\subseteq\bigcup_{i=1}^{t}E_{nn}(\alpha_{s_{i}})$
thus
\begin{equation*}
 S_{k,s_{1}}\wedge\dots\wedge S_{k,s_{t}}=(\mu_{k})_{*}(\beta_{k,s_{1}}\wedge\dots\wedge\beta_{k,s_{t}})~\text{on}~X.
\end{equation*}
So we can write the last equality as follows
\begin{equation}\label{eq1}
 T_{k,s_{1}}\wedge\dots\wedge T_{k,s_{t}} + O(\varepsilon_k)=(\mu_{k})_{*}(\beta_{k,s_{1}}\wedge\dots\wedge\beta_{k,s_{t}})
\end{equation}
thanks to uniform control of the mass by cohomology classes.
 Passing in cohomology in (\ref{eq1}) and letting  $k$  to $+\infty$, the statement holds.\hfill$\square$ 
\\
\flushleft Now we want to study the extremal case of Proposition (\ref{pro2}), i.e. one can consider the difference 
$\alpha_{1}\cdots\alpha_{q}-\langle\alpha_{1}\cdots\alpha_{q}\rangle$. The case $q=1$ is given by \cite{Bou1}.
We define for every $0\leq q\leq n$
\begin{equation*}
 \mathcal{M}^{q}:=\{\alpha\in\mathcal{E}:~\dim(E_{nn}(\alpha))\leq n-q\},
\end{equation*}
then we have the following
\begin{prop}
 $\mathcal{M}^{q}$ is a convex closed cone of $\mathcal{E}$ for all $q=0,\dots,n$.
\end{prop}
\textit{Proof.} $\mathcal{M}^{q}$ is convex as it follows from the convexity of the map $\mathcal{E}\to\mathbb{R}$, $\alpha\to\nu(\alpha,x)$.\
Now let $\alpha$ be in the closure of $\mathcal{M}^{q}$, for $\varepsilon>0$ small enough $\alpha+\varepsilon\omega$ is in the interior of $\mathcal{M}^{q}$.
By using Remark (\ref{rmk9}) one has that 
\begin{equation*}
 \displaystyle E_{nn}(\alpha)=\bigcup_{\varepsilon>0,\varepsilon\in\mathbb{Q}}E_{nn}(\alpha+\varepsilon\omega)
\end{equation*}
thus $\dim (E_{nn}(\alpha))\leq n-q$ then $\mathcal{M}^{q}$ is closed.\\
Finally we observe that for $q=0,1$ $\mathcal{M}^{q}=\mathcal{E}$, while for $q=n$ $\mathcal{M}^{n}=\mathcal{N}$ the cone of nef classes.
\hfill$\square$ \\
\
\newline
As a consequence of the definition of non-pluripolar product one has
\begin{cor}\label{cor1}
 If $T_{1},\dots,T_{q}$ are closed positive $(1,1)-currents$ such that the union $F$ of their unbounded locus is contained in an analytic set $Y$ of pure dimension $n-q$.
Then
\begin{equation}
 \nu(\langle T_{1}\wedge\cdots\wedge T_{q}\rangle,Z)=0
\end{equation}
for all irreducible components $Z$ of dimension $n-q$ in $Y$.
\end{cor}
\textit{Proof.} Let's consider the Siu decomposition of $\langle T_{1}\wedge\cdots\wedge T_{q}\rangle=R+\sum_{h}\nu(\langle T_{1}\wedge\cdots\wedge T_{q}\rangle,Z_{h})[Z_h]$
then $Z_{h}\subseteq Y$ and we find that $\chi_{X-Y}\langle T_{1}\wedge\cdots\wedge T_{q}\rangle=\chi_{X-Y}R$, then the extension by zero of 
$\chi_{X-Y}\langle T_{1}\wedge\cdots\wedge T_{q}\rangle$ on $Y$ is $\langle T_{1}\wedge\cdots\wedge T_{q}\rangle$ then $R\geq\langle T_{1}\wedge\cdots\wedge T_{q}\rangle$ 
and since $\nu(R,Z)=0$ for every analytic set of dimension $n-q$ we infer the statement.\hfill$\square$
\\
\begin{prop}\label{pro9}
 Let $T_{j},T_{j}^{'}$ be positive closed currents for $j=1,\dots,q$ with $T_{j}\preceq T_{j}^{'}$, such that
$T_{1}\wedge\cdots\wedge T_{q}$ and $T_{1}^{'}\wedge\cdots\wedge T_{q}^{'}$ are defined. Then
\begin{equation}\label{eq13}
 \nu(T_{1}\wedge\cdots\wedge T_{q},x)\leq\nu(T_{1}^{'}\wedge\cdots\wedge T_{q}^{'},x)
\end{equation}
 for all $x\in X.$\\
In particular if $T_{min,j}$ and $T_{min,j}^{'}$ are currents with minimal singularities in $\alpha_{j}$ then 
$\nu(T_{min,1}\wedge\cdots\wedge T_{min,q},x)$ does not depend in the chosen currents with minimal singularities.
\end{prop}
\textit{Proof.} 
 We can choose local coordinates $z=(z_1,\cdots,z_n)$ centered in $x$.
Let $v_{j}\leq0$ and $u_{j}\leq0$ be the local potentials for $T_{j}$ and $T_{j}^{'}$ respectively, then $v_{j}+C_{j}\geq u_{j}$ for
some constants $C_j\geq0$.  Let's define
\begin{itemize}
\item $ v_{j}^{'}:=v_{j}+\lambda\log|z|$\\
\item $u_{j}^{'}:=u_{j}+\lambda\log|z|$
\end{itemize}
for $\lambda>0$.  Thus we find
\begin{equation*}
 v_{j}^{'}+C_{j}\geq u_{j}^{'}
\end{equation*}
then dividing by $u_{j}^{'}$ and passing to the $\limsup$ we get
\begin{equation}
 \limsup_{z\to x}\dfrac{v_{j}^{'}(z)}{u_{j}^{'}(z)}\leq1.
\end{equation}

Applying (\cite{Dem2}, Theorem 5.9) we find
\begin{equation}\label{eq12}
 \displaystyle\nu(\bigwedge_{j=1}^{q}dd^{c}v_{j}^{'},x)\leq\nu(\bigwedge_{j=1}^{q}dd^{c}u_{j}^{'},x).
\end{equation}
 In the calculation of the left-side hand in (\ref{eq12}) we get a polynomial  $p(\lambda)$ whose constant term $a_{0}$ is just $q$-times the wedge product
of $dd^{c}v_{j}$, while for the right-side hand we get a polynomial $q(\lambda)$ whose constant term $b_{0}$ is just $q$-times the wedge product
of $dd^{c}u_{j}$. Then letting $\lambda\to0$ equation (\ref{eq13})
holds.\hfill $\square$

\begin{rmk}\label{rmk8}
 If $\alpha_{1},\dots,\alpha_{q}\in\mathcal{M}^{q}$ with  then  $\beta_{i}=\alpha_{i}+\varepsilon\omega$ are in the $\text{int}(\mathcal{M}^{q})$ for all $1\leq i\leq q$.\\
If $0<\varepsilon^{'}<\varepsilon$ then $T_{min,\alpha_{j}+\varepsilon\omega}\preceq T_{min,\alpha_{j}+\varepsilon^{'}\omega}+(\varepsilon-\varepsilon^{'})\omega$
then by Proposition (\ref{pro9})  $\nu(\bigwedge_{j=1}^{q}T_{min,\alpha_{j}+\varepsilon\omega},x)$ increases when $\varepsilon\to0$
for all $x\in X$ and 
it is bounded in terms of $\alpha_{1}\cdots\alpha_{q}\cdot\{\omega\}.$\end{rmk}
Then we can give the following
\begin{defn}\label{defini}
 The ``multiplicity'' of $\alpha_{1},\dots,\alpha_{q}$ at $x\in X$ is defined as
\begin{equation}
 \displaystyle\nu(\alpha_{1}\cdots\alpha_{q},x):=\sup_{\varepsilon>0}\nu(\bigwedge_{j=1}^{q}T_{min,\alpha_{j}+\varepsilon\omega},x)
\end{equation}
and we define $\nu(\alpha_{1}\cdots\alpha_{q}, Y)=\inf\{\nu(\alpha_{1}\cdots\alpha_{q},x),~x\in Y\}$.
\end{defn}

\begin{prop}\label{equali}
 If $\alpha_{1},\dots,\alpha_{q}$ are in the interior of $\mathcal{M}^{q}$.
Then
\begin{equation*}
 \nu(\alpha_{1}\cdots\alpha_{q},x)=\nu(T_{min,\alpha_{1}}\wedge\cdots\wedge T_{min,\alpha_{q}},x),
\end{equation*}
where $T_{min,\alpha_{i}}$ is a positive current with minimal singularities in $\alpha_{i}.$
\end{prop}
\textit{Proof.} It is clear that $\nu(\alpha_{1}\cdots\alpha_{q},x)\leq\nu(T_{min,\alpha_{1}}\wedge\cdots\wedge T_{min,\alpha_{q}},x).$
For the other inequality let $S_{j}\in\alpha_{j}$ be a K$\ddot{\text{a}}$hler current such that $S_{j}\geq\omega$ and let $T_{min,\varepsilon,j}$ be
a minimal  current in $\alpha_{j}[-\varepsilon\omega]$. Then by remark (\ref{rmk9}) we have that
\begin{equation*}
\nu(\lambda S_{j}+(1-\lambda)T_{min,\varepsilon,j},x)=\nu(\lambda S_{j}+(1-\lambda)(T_{min,\alpha_{j}+\varepsilon\omega}-\varepsilon\omega),x)
\end{equation*}
where the current $K_{j}=\lambda S_{j}+(1-\lambda)(T_{min,\alpha_{j}+\varepsilon\omega}-\varepsilon\omega)\geq0$ is a positive current representing $\alpha_{j}$, with
$\lambda=\frac{\varepsilon}{1+\varepsilon}$. 
By using Proposition (\ref{pro9}) we infer that
\begin{equation*}
 \nu(K_{1}\wedge\cdots\wedge K_{q},x)\geq\nu(T_{min,\alpha_{1}}\wedge\cdots\wedge T_{min,\alpha_{q}},x).
\end{equation*}
Now by using the uniform bound of mass in terms of cohomology classes and letting $\varepsilon\to0$ (hence $\lambda\to0$), we obtain
the other inequality.

\hfill$\square$

Before proving \textbf{Theorem 1}, we recall that 
\begin{lem}\label{inter}
 If $\alpha_{1},\dots,\alpha_{q}$ are in $\mathcal{M}^{q}$ then
\begin{equation}\label{inc1}
\displaystyle \bigcup_{\varepsilon>0,\varepsilon\in\mathbb{Q}}\bigcap_{j=1}^{q}E_{nK}(\alpha_{j}+\varepsilon\omega)=\bigcap_{j=1}^{q}\bigcup_{\varepsilon>0,\varepsilon\in\mathbb{Q}}E_{nK}(\alpha_{j}+\varepsilon\omega)=
\bigcap_{j=1}^{q}E_{nn}(\alpha_{j}).
\end{equation}
\end{lem}
\textit{Proof.} For $\varepsilon>0$ $\alpha_{j}+\varepsilon\omega\in\text{int}(\mathcal{M}^{q})$ (interior of $\mathcal{M}^{q}$) for all
$j=1,\dots,q$, in particular $\alpha_{j}+\frac{\varepsilon}{2}\omega\in\text{int}(\mathcal{M}^{q}).$ Let $\{T_{k,j}\}\in\alpha_{j}+\frac{\varepsilon}{2}\omega$
be a sequence of currents with analytic singularities as in Theorem (\ref{thm2}) such that $T_{k,j}\rightharpoonup T_{min,\alpha_{j}+\frac{\varepsilon}{2}\omega}$
where $T_{min,\alpha_{j}+\frac{\varepsilon}{2}\omega}$ is a positive current with minimal singularities in $\alpha_{j}+\frac{\varepsilon}{2}\omega.$\\
Now for $k>>0$ $T_{k,j}+\frac{\varepsilon}{2}\omega$ are K$\ddot{\text{a}}$hler currents with analytic singularities, representing (as currents)
$\alpha_{j}+\varepsilon\omega,$ then we have the following inclusions
\begin{equation}\label{inc2}
 E_{nn}(\alpha_{j}+\varepsilon\omega)\subseteq E_{nK}(\alpha_{j}+\varepsilon\omega)\subseteq E_{+}(T_{k,j}+\frac{\varepsilon}{2}\omega)\subseteq E_{nn}(\alpha_{j}+\frac{\varepsilon}{2}\omega)
\end{equation}
thus it follows 
\begin{equation}\label{inc3}
 \displaystyle E_{nn}(\alpha_{j})=\bigcup_{\varepsilon>0,\varepsilon\in\mathbb{Q}}E_{nK}(\alpha_{j}+\varepsilon\omega).
\end{equation}
Now thanks to (\ref{inc3}) the equality (\ref{inc1}) is equivalent to the following
\begin{equation}\label{inc4}
 \displaystyle\bigcup_{\varepsilon>0,\varepsilon\in\mathbb{Q}}\bigcap_{j=1}^{q}E_{nK}(\alpha_{j}+\varepsilon\omega)=
\bigcap_{j=1}^{q}\bigcup_{\varepsilon>0,\varepsilon\in\mathbb{Q}}E_{nn}(\alpha_{j}+\varepsilon\omega).
\end{equation}
Let's define $A$ and $B$ the left-hand side and the right-hand side of (\ref{inc4}).\\
Now if $x\in A\Rightarrow\exists\varepsilon_{0}>0$ such that $\forall j=1,\dots,q,$ $x\in\bigcap_{j=1}^{q}E_{nK}(\alpha_{j}+\varepsilon_{0}\omega)$
but $E_{nK}(\alpha_{j}+\varepsilon_{0}\omega)\subseteq E_{nn}(\alpha_{j}+\frac{\varepsilon_{0}}{2}\omega)$ $\forall j=1,\dots,q$ thanks to (\ref{inc2}).
Then $x\in\bigcup_{\varepsilon>0}E_{nn}(\alpha_{j}+\varepsilon\omega)$ $\forall j=1,\dots,q$, i.e. $x\in B.$ \\
Viceversa if $x\in B$ then there exist $\varepsilon_{1},\dots,\varepsilon_{q}>0$ such that $x\in E_{nn}(\alpha_{j}+\varepsilon_{j}\omega)$
for $j=1,\dots,q$. Let $\bar{\varepsilon}=\min_{j}\varepsilon_{j}$ then again thanks to (\ref{inc2}) we have
\begin{equation*}
 E_{nn}(\alpha_{j}+\varepsilon_{j}\omega)\subseteq E_{nn}(\alpha_{j}+\bar{\varepsilon}\omega)\subseteq E_{nK}(\alpha_{j}+\bar{\varepsilon}\omega).
\end{equation*}
It follows that $x\in E_{nk}(\alpha_{j}+\bar{\varepsilon}\omega)$ for all $j=1,\dots,q$ then $x\in A.$ \hfill$\square$
\\
\
\newline
Using the previous Lemma we give a proof of the following well-known
\begin{prop}\label{zero}
 The set $E_{nn}(\alpha)$ does not have irreducible components of dimension zero if $\alpha\in\mathcal{E}$.
\end{prop}
\textit{Proof.} Let $x\in E_{nn}(\alpha)$ be an irreducible zero dimensional component. By Lemma (\ref{inter}) we have
\begin{equation*}
 E_{nn}(\alpha)=\bigcup_{\varepsilon>0,\varepsilon\in\mathbb{Q}}E_{nK}(\alpha+\varepsilon\omega)
\end{equation*}
then $x$ is an isolated point in $E_{nK}(\alpha+\varepsilon\omega)$ for all $\varepsilon>0$ rational since $E_{nK}(\alpha+\varepsilon\omega)$
are analytic sets. Now let $\{T_{k}\}\in\alpha$ be a sequence of currents with analytic singularities as in Theorem (\ref{thm2}) such
that $T_{k}\rightharpoonup T_{min}$ weakly where $T_{min}$ is a positive current with minimal singularities in $\alpha.$ Then for $k>>0$
$T_{k}+\varepsilon\omega$ are K$\ddot{\text{a}}$hler currents with analytic singularities only in $x$ representing $\alpha+\varepsilon\omega$ (as currents)
; hence by the gluing property of plurisubharmonic functions we can find new K$\ddot{\text{a}}$hler currents $S_{k}+\varepsilon\omega$ still representing $\alpha+\varepsilon\omega$ 
as currents, but with $\nu(S_{k}+\varepsilon\omega,x)=0$, a contradiction. \hfill$\square$

\newpage
Now we can prove
\begin{thm}\label{pro3}
 Let $X$ be a compact K$\ddot{\text{a}}$hler manifold and be $\alpha_{1},\dots,\alpha_{q}\in{\mathcal{M}}^{q}$.
Let $\{Y_{h}\}_{h\in H}$ the family of codimension $q$ components of the set $\bigcap_{j=1}^{q}E_{nn}(\alpha_{j})$.
Then we have the following decomposition:\\

\item i) $\displaystyle\alpha_{1}\cdots\alpha_{q}=\langle\alpha_{1}\cdots\alpha_{q}\rangle+\left\{\sum_{h\in H}\nu(\alpha_{1}\cdots\alpha_{q},Y_{h})[Y_{h}]\right\}$,
 with \hspace*{2mm} $\nu(\alpha_{1}\cdots\alpha_{q},Y_{h})>0$ $\forall h\in H;$ 
\\

\item ii) if $\alpha_{j}\in \text{int}(\mathcal{M}^{q})$ for $1\leq j\leq q$, then the set $H$ is finite and \hspace*{2mm}
$\nu(\alpha_{1}\cdots\alpha_{q},Y_{h})=\nu(T_{min,\alpha_{1}}\wedge\cdots\wedge T_{min,\alpha_{q}},Y_{h})$ for all $h\in H$.
\end{thm}
\textit{Proof.} $ii)$ By the proof of Lemma (\ref{inter}) for $\varepsilon>0$ small enough $\alpha_{j}-\varepsilon\omega\in\text{int}(\mathcal{M}^{q})$ and
let $\{T_{k,j}\}\in\alpha_{j}-\varepsilon\omega$ be a sequence of currents with  analytic singularities as in Theorem (\ref{thm2}) such that $T_{k,j}\rightharpoonup T_{min,\alpha_{j}-\varepsilon\omega}$
then 
\begin{equation*}
 \nu(T_{k,j}+\varepsilon\omega,x)=\nu(T_{k,j},x)\leq\nu(T_{min,\alpha_{j}-\varepsilon\omega},x)
\end{equation*}
But $\alpha_{j}-\varepsilon\omega$ is still in $\text{int}(\mathcal{M}^{q})$ thus $E_{+}(T_{k,j}+\varepsilon\omega)$ is an analytic set of dimension $\leq n-q$
so the set  of components of dimension $n-q$ of $E_{nn}(\alpha_j)$ and of $E_{nK}(\alpha_j)$ is finite
. Let $T_{min,j}$ be a positive current with minimal singularities in $\alpha_j$, since $\displaystyle\dim_{\mathcal{H}}\left(\bigcup_{i}^{q}E_{nK}(\alpha_{j})\right)\leq2n-2q$
($\dim_{\mathcal{H}}$ is the Hausdorff dimension), then $T:=T_{min,1}\wedge\dots\wedge T_{min,q}$ is well defined on $X$ according to (\cite{Dem1}) and
by Proposition (\ref{equali}) $\nu(\alpha_{1}\cdots\alpha_{q},x)=\nu(T_{min,1}\wedge\cdots\wedge T_{min,q},x)$ for all $x\in X$.
Let's consider the positive current $T-\langle T\rangle:$ this difference is identically zero outside $\displaystyle\bigcup_{j}^{q}E_{nK}(\alpha_{j})$ and  by the support theorem of positive
currents it follows:
\begin{equation*}
 \displaystyle T-\langle T\rangle=\sum_{s=1}^{t}\lambda_{s}[Y_s]
\end{equation*}
where $Y_s$ are the irreducible components of codimension $q$ of $\displaystyle\bigcup_{j}^{q}E_{nK}(\alpha_{j})$ and $\lambda_s\geq0$.
Now by Corollary (\ref{cor1})
\begin{equation}\label{eq14}
 \displaystyle T=\langle T\rangle+\sum_{h=1}^{t}\nu(T,Y_h)[Y_h]
\end{equation}
is the Siu decomposition of $T$. Now if $Y_{s_{0}}$ is  contained in $\bigcap_{j=1}^{q}E_{nn}(\alpha_{j})$ then by \cite{Dem2}
\begin{equation*}
 \displaystyle\nu(T_{min,1}\wedge\cdots\wedge T_{min,q},Y_{s_{0}})\geq\prod_{j=1}^{q}\nu(T_{min,j},Y_{s_{0}})>0.
\end{equation*}
On the other hand let $Y_{s}\subseteq\bigcup_{j=1}^{q}E_{nK}(\alpha_{j})$ be an irreducible component of codimension $q$ not entirely
contained in $\bigcap_{j=1}^{q}E_{nn}(\alpha_{j})$. Then by Lemma (\ref{inter}) there exists $j_{0}$ such that for all $\varepsilon>0$
$Y_{s}$ is not entirely contained in $\bigcup_{\varepsilon>0,\varepsilon\in\mathbb{Q}}E_{nK}(\alpha_{j_{0}}+\varepsilon\omega)$.
It is not restrictive to suppose $j_{0}=1$ then the local potential $\varphi_{min,1}$ of $T_{min,\alpha_{1}+\varepsilon\omega}$ is locally
bounded in a neighborhood of a generic point $x\in Y_{s}$.\\
Let $B=B(x,r)$ be an open ball of center $x$ and radius $r>0$ small enough such that $B$ is contained in a coordinate patch with local
coordinates $z=(z_{1},\dots,z_{n})$; set $\beta=dd^{c}|z|^{2}$ and $S=\bigwedge_{i=2}^{q}dd^{c}\varphi_{min,i}$. Let $\chi$ be a smooth
function with compact support in $B$ such that $0\leq\chi\leq1$. Then
\begin{equation*}
 \int_{B}\chi dd^{c}\varphi_{min,1}\wedge S\wedge\beta^{n-q}=\int_{B}\varphi_{min,1}\cdot S\wedge dd^{c}\chi\wedge\beta^{n-q}
\end{equation*}
where the equality is obtained integrating by parts and using that $S$ and $\beta$ are closed.\\
 By using Chern-Levine-Nirenberg inequalities (\cite{Dem2}, Proposition 1.3),  we get
\begin{equation*}
 \int_{B}\varphi_{min,1}\cdot S\wedge dd^{c}\chi\wedge\beta^{n-q}\leq \text{Vol}(B)||S||_{X}C
\end{equation*}
where $C$ is the product of the bound for the coefficients of the smooth form $dd^{c}\chi\wedge\beta^{n-q}$ and of the bound for $\varphi_{min,1}$,
Vol$(B)$ is the volume of $B$ and $||S||_{X}$ is the mass of the current $S$ that can be bounded in terms of $\alpha_{2}\cdots\alpha_{q}\cdot\{\omega\}$.
Using the fact that Vol$(B)\sim r^{2n}$ and that
\begin{equation*}
 \displaystyle\nu(dd^{c}\varphi_{min,1}\wedge S,x)=\lim_{r\to0}\frac{1}{2^{2(n-q)}r^{2(n-q)}}\int_{B(x,r)}dd^{c}\varphi_{min,1}\wedge S\wedge\beta^{n-q}
\end{equation*}
we can infer that $\nu(T,Y_{s})=0.$
 Then  in equation (\ref{eq14}) $Y_{h}$ are actually irreducible components of codimension $q$ of 
$\bigcap_{j=1}^{q}E_{nn}(\alpha_{j})$. 
 
Then passing in cohomology and letting $\varepsilon\to0$ the statement holds. \\
To prove $i)$ it is sufficient to observe that $\alpha_{j}+\varepsilon\omega\in\mathring{\mathcal{M}}^{q}$ for $1\leq j\leq q$,
applying $ii)$ to $\alpha_{j}+\varepsilon\omega$  take the limit as $\varepsilon\to0$ and using Lemma (\ref{inter}).
\hfill $\square$\\

\flushleft Let's observe that if one of $\alpha_i$'s is big but it is on the boundary of $\mathcal{M}^{q},$ it may be not true that the number of irreducible
components of codimension $q$ of $E_{nn}(\alpha_i)$ is finite. Since it could happen that $E_{nK}(\alpha_{i})$ have irreducible
components of codimension $q-1$, hence an infinite number of components of codimension $q$ of $E_{nn}(\alpha_i)$ could be contained in one of these
components.\\
An immediate consequence of Theorem (\ref{pro3}) is the following
\begin{cor}
 If $\alpha\in\emph{int}(\mathcal{M}^{q})$ then $\nu(T_{min}^{q},Y)=0$ for every irreducible analytic set of codimension $q$ not completely
contained in $E_{nn}(\alpha)$.
\end{cor}

Another consequence of Theorem (\ref{pro3}) is a partial  converse of Theorem (\ref{thm3} (ii))
\begin{prop}
 Let $(X,\omega)$ be a compact K$\ddot{\text{a}}$hler manifold of complex dimension $n$.Then:
\begin{equation*}
\alpha\in\mathcal{M}^{q}~\text{if and only if}~\alpha^{j}=(\alpha^{j})~\text{for all}~1\leq j\leq q-1.
\end{equation*}
In particular, $\alpha$ is nef if and only if $\alpha^{j}=(\alpha^{j})$ for all $1\leq j\leq n-1.$
\end{prop}
\textit{Proof.} If $\alpha\in\mathcal{M}^{q}$  by Proposition
(\ref{pro2}) $\alpha^j=(\alpha^{j})$ for all $1\leq j<q$.\\
To prove the vice versa we go by induction on $q$.\\
For $q=2$, the statement follows from \cite{Bou1}. Assume now the statement for a given $q$ and let $\alpha^{j}=(\alpha^{j})$ for
$1\leq j\leq q+1$, then by induction we know that $\alpha\in\mathcal{M}^{q+1},$ the by Theorem (\ref{pro3}) we know that
cod$(E_{nn}(\alpha))\geq q+2$, i.e. $\alpha\in\mathcal{M}^{q+2}.$\\
The last statement follows from Proposition (\ref{zero}). \hfill$\square$
\\
\section{Algebraic Morse inequalities in codimension ``\textit{s}''}

Let $X$ be a compact K$\ddot{\text{a}}$hler manifold, $E$ a holomorphic vector bundle of rank $r$ and $L$ a line bundle over $X$. If
$L$ is equipped with a smooth metric of curvature form $\Theta(L)$, one defines the $s-$index set of $L$ to be the open subset
\begin{equation*}
 \displaystyle X(s,L)=\left\{x\in X~:~i\Theta(L)_{x}~\text{has}\begin{array}{cl}
                                                                s &\text{negative eigenvalues}\\
                                                                n-s &\text{positive eigenvalues}
                                                               \end{array}\right\}
\end{equation*}
for $0\leq s\leq n=\dim_{\mathbb{C}}(X)$. Hence $X$ admits a partition $X=\Delta\cup\bigcup_{s}X(s,L)$ where $\Delta=\{x\in X~:~\det(\Theta(L)_{x})=0\}.$
One also introduces
\begin{equation*}
\displaystyle X(\leq s,L)=\bigcup_{0\leq j\leq s}X(j,L).
\end{equation*}
It is shown in \cite{Dem5} that the cohomology groups $H^{s}(X,E\otimes\mathcal{O}(kL))$ satisfy the following asymptotic \textit{weak Morse inequalities}
as $k\to+\infty$
\begin{equation}
\displaystyle h^{s}(X,E\otimes\mathcal{O}(kL))\leq r\frac{k^{n}}{n!}\int_{X(s,L)}(-1)^{s}\left(\frac{i}{2\pi}\Theta(L)\right)^{n}+o(k^n).
\end{equation}
A sharper form is given by the \textit{strong Morse inequalities}
\begin{equation}\label{strongM}
 \displaystyle\sum_{j=0}^{s}(-1)^{s-j}h^{j}(X,E\otimes\mathcal{O}(kL))\leq r\frac{k^{n}}{n!}\int_{X(\leq s,L)}(-1)^{s}\left(\frac{i}{2\pi}\Theta(L)\right)^{n}+o(k^n).
\end{equation}
One difficulty in the application of these inequalities is that the curvature integral is in general quite uneasy to compute, it is neither a topological
nor an algebraic invariant. However a special case of the Morse inequalities can be reformulated in a more algebraic setting in which only algebraic invariants
are involved, se e.g. \cite{Trap1}, \cite{Siu2}.
\newpage
Here we give an \emph{algebraic} reformulation for the Holomorphic Morse inequalities proved by J.-P. Demailly (\cite{Dem5}) in a more
general setting.

\begin{thm}[\cite{Dem6}]
 Let $V=L-F$ be a holomorphic line bundle over a compact K$\ddot{\text{a}}$hler manifold $X$, where $L$ and $F$ are nef line bundle.
Then for every $s=0,\dots,n$, there is an asymptotic strong Morse inequality
\begin{equation*}
 \displaystyle\sum_{0\leq j\leq s}(-1)^{s-j}h^{j}(X,kV)\leq\frac{k^n}{n!}\sum_{0\leq j\leq s}(-1)^{s-j}\binom{n}{j}L^{n-j}F^{j}+o(k^n).
\end{equation*}
\end{thm}

If $F$ is not nef but it is just pseudo-effective with codim$(E_{nn}(F))\leq s+1$, then
\begin{thm}[\cite{Trap2}]
 Let $L$ and $F$ be holomorphic line bundle over $X$ a projective compact manifold, with $L$ nef and $F\in\mathcal{E}$ such that $\dim(E_{nn}(F))\leq n-(s+1)$. Then
 for $0\leq m\leq s$ we have the following holomorphic Morse inequalities
\begin{equation}
 \displaystyle\sum_{0\leq j\leq m}(-1)^{m-j}h^{j}(X,k(L-F))\leq\frac{k^n}{n!}\sum_{0\leq j\leq m}(-1)^{m-j}\binom{n}{j}L^{n-j}F^{j}+o(k^n)
\end{equation}
\end{thm}

However, when $s=1$ there is a version of the algebraic Morse inequalities which uses the full divisorial Zariski decomposition of $F$
\begin{thm}[\cite{Trap2}]\label{thm5}
 Let $L$ and $F$ be line bundles over $X$, assume $L$ nef and $F$ pseudo-effective, let $F=(F)+\{N(F)\}$ be the divisorial Zariski decomposition of $F$
with $N(\{F\})=\sum_{j=1}^{N}\nu(F,D_{j})[D_{j}]$, and let $\{u\}$ be a nef cohomology class in $H^{2}(X,\mathbb{R})$ such that
$c_{1}(\mathcal{O}_{{T_X}}(1))+\pi^{*}\{u\}$ is a nef cohomology class in $H^{2}(\mathbb{P}(T_{X}^{*}),\mathbb{R})$. Then
\begin{equation}
 \displaystyle \limsup_{k\to+\infty}\frac{n!}{k^{n}}h^{0}(X,k(L-F))\geq L^{n}-nL^{n-1}(F)-n\sum_{j=1}^{N}(L+\nu(F,D_{j})\{u\})^{n-1}\nu(F,D_{j})[D_j]
\end{equation}
 \end{thm}

We can also treat the extremal case, i.e. when the codimension of $E_{nn}(F)=s$. We give two formulations.
\begin{thm}[First Formulation]\label{Mor1}
 Let $X$ be a compact projective manifold of complex dimension $n$. Let $L$ and $F$ be two line bundles over $X$ with $L$ nef and
$F\in\mathcal{M}^{s}$ with $\dim (E_{nn}(F))=n-s$. Let $\{Y_{t}\}_{t\in T}$ be the irreducible components (possibly infinite) of codimension s of 
$E_{nn}(F)$, and
let $\nu_{t}$ and $\nu_{t}^{'}$  be the multiplicities of $F^{s}$ and the multiplicity of $F$ along
$Y_{t}$ respectively. Then we have the following Morse inequalities
\begin{multline}
 \displaystyle \sum_{j=0}^{s}(-1)^{s-j}h^{j}(X,k(L-F))\leq
 \frac{k^{n}}{n!}(\sum_{j=0}^{s}\binom{n}{j}(-1)^{s-j}L^{n-j}(F^{j})+\\
 \sum_{t=1}^{+\infty}\binom{n}{s}(L+\nu_{t}^{'}\{u\})^{n-s}\nu_{t}[Y_{t}])+o(k^n).
\end{multline}
\end{thm}
\textit{Proof.}
We first assume $L$ to be ample and $F\in\text{int}(\mathcal{M}^{s})$  .
Since $F\in\text{int}(\mathcal{M}^{s})$, in particular $F$ is big,  then if $T_{min}\in c_1(F)$ is a current with minimal singularities then  we have
$\nu^{'}(c_1(F),x)=\nu^{'}(T_{min},x)$ for all $x\in X$.
We recall that there is a finite number of $Y_{j}$'s  thanks to Theorem (\ref{pro3}).
Set $0<\nu_{1}^{'}<\cdots<\nu_{h}^{'}$, where $\nu_{1}^{'}$ is the minimum of $\nu_{j}^{'}$ positive and so on. Let's choose $c_{i}$ real positive numbers such that
\begin{equation}
 0<c_{0}<\nu_{1}^{'}<c_{1}<\nu_{2}^{'}<\cdots<\nu_{h-1}^{'}<c_{h-1}<\nu_{h}^{'}<c_{h}<\nu:=\max_{x\in X}(\nu(T_{min},x))+1.
\end{equation}
  Choose a positive closed $(1,1)-$form $\omega$ which is the curvature  of a smooth metric
on the ample line bundle $L$. Now by Theorem (\ref{thm2}) we know that there exist a sequence of closed smooth forms $T_k$ in the cohomology
class of $c_1(F)$, a decreasing sequence of positive functions $\lambda_k(x)$, and a decreasing sequence of positive
real numbers $\varepsilon_k$, with the following properties
\begin{itemize}
\item $T_k\rightharpoonup T_{min}$,\\
\item $\lambda_{k}(x)$ converges to $\nu(T_{min},x)$ for all $x\in X$,\\
\item $\varepsilon_{k}$ converges to $0$,\\
\item $T_{k}>-\lambda_{k}u-2\varepsilon_{k}\omega.$
\end{itemize}
Let $c$ be a positive real number, then one can defines $\Omega_{k,c}=\{x\in X~:~\lambda_{k}(x)<c\}$, $v_{k,c}=2\varepsilon_{k}+cu$ and
$w_{k,c}=v_{k,c}+\omega.$ Then on $\Omega_{k,c}$ the forms $T_k+v_{k,c}$ and $w_{k,c}$ are positive. Now $\omega-T_k$ is the curvature of smooth
metric on $L-F$. Let $\alpha_{1}\leq\cdots\alpha_{n}$ be the eigenvalues of $\omega-T_k$ with respect to $w_{k,c}$, so that $\alpha_{j}\leq1$ for all $j.$
Let $X(p)$ be the $p$-index set for $L-F$ and $X(\leq p)$ be the set of points of index at most $p$ . Then on $X(p)$ we have $(-1)^{p}(\omega-T_{k})^{n}\leq(-1)^{p}\alpha_{1}\cdots\alpha_{n}w_{k,c}^{n}$.
Now
\begin{equation*}
\binom{n}{p}w_{k,c}^{n-p}\wedge(T_{k}+v_{k,c})^{p}=\binom{n}{p}w_{k}^{n-p}\wedge(w_{k,c}-(\omega-T_k))^{p}=\sigma_{p}(1-\alpha)w_{k,c}^{n},
\end{equation*}
where $\sigma_{p}(1-\alpha)$ is the $p$-th elementary symmetric function in $1-\alpha_{1},\dots,1-\alpha_{n}$.
However, since $\alpha_j<1$ for $1\leq j\leq n$, it follows that on $X(p)$ we have
\begin{equation*}
 \sigma_{p}(1-\alpha)\geq(1-\alpha_{1})\cdots(1-\alpha_{p})\geq(-1)^{p}\alpha_{1}\cdots\alpha_{p}.
\end{equation*}
Furthermore one can easly prove by induction on $n$ (\cite{Dem6}) that
\begin{equation*}
\displaystyle\sum_{j=0}^{p}(-1)^{p-j}\sigma_{p}(1-\alpha)\leq\chi_{X(\leq p)}(-1)^{p}\alpha_{1}\cdots\alpha_{n},
\end{equation*}
where $\chi_{X(\leq p)}$ is the characteristic function of the set $X(\leq p)$.

Now we want to estimate $\int_{X(\leq s)}(-1)^{s}(\omega-T_k)^{n}$.
\begin{equation}\label{eq3}
 \int_{X(\leq s)}(-1)^{s}(\omega-T_k)^{n}=\int_{X(\leq s)\cap\Omega_{k,c_{0}}}(-1)^{s}(\omega-T_k)^{n}+\int_{X(\leq s)\cap\Omega_{k,c_{0}}^{c}}(-1)^{s}(\omega-T_k)^{n}.
\end{equation}\newpage
$\boxed{\text{Estimate of}~ \int_{X(\leq s)\cap\Omega_{k,c_{0}}}(-1)^{s}(\omega-T_k)^{n}}$\\
\flushleft We get
\begin{equation}
 \displaystyle\int_{X(\leq s)\cap\Omega_{k,c_{0}}}(-1)^{s}(\omega-T_k)^{n}\leq\int_{X(\leq s)\cap\Omega_{k,c_{0}}}\sum_{j=0}^{s}\binom{n}{j}(-1)^{s-j}w_{k,c_{0}}^{n-j}\wedge(T_k+v_{k,c_{0}})^{j},
\end{equation}
Set $S_k:=\chi_{\Omega_{k,c_{0}}}w_{k,c_{0}}^{n-j}\wedge(T_k+v_{k,c_{0}})^{j}$ then
\begin{equation}
 S_k\rightharpoonup\chi_{E_{c_{0}}(T_{min})^{c}}(\omega+c_{0}u)^{n-j}\wedge(T_{min}+c_{0}u)^{j}
\end{equation}
 where $E_{c_{0}}(T_{min})^{c}$ is the complementary set of $E_{c_{0}}(T_{min})$.\\
Now for $j=0,\dots,s-1$, since cod$(E_{c_{0}}(T_{min}))=s$ and $T_{min}^{j}=\langle T_{min}^{j}\rangle$ by Proposition (\ref{pro2}), then
\begin{multline}\label{mult1}
 \displaystyle\int_{X}\chi_{E_{c_{0}}(T_{min})^{c}}(\omega+c_{0}u)^{n-j}\wedge(T_{min}+c_{0}u)^{j}=\\
\int_{X}\chi_{E_{c_{0}}(T_{min})^{c}}(\omega+c_{0}u)^{n-j}\wedge\left(\sum_{h=0}^{j}\binom{j}{h}\langle T_{min}^{h}\rangle\wedge(c_{0}u)^{j-h}\right)=\\
\int_{X}(\omega+c_{0}u)^{n-j}\wedge\left(\sum_{h=0}^{j}\binom{j}{h}\langle T_{min}^{h}\rangle\wedge(c_{0}u)^{j-h}\right).
\end{multline}
For $j=s$, by Theorem (\ref{pro3}) one has $T_{min}^{s}=\langle T_{min}^{s}\rangle+\sum_{t=1}^{N}\nu(T_{min}^{s},Y_t)[Y_t]$, where
$Y_t$ are the irreducible components of codimension $s$ in $E_{nn}(T_{min})$ . Then
we obtain the same identity for $j=s$ as in (\ref{mult1}). 
So first letting $k\to+\infty$ and then $c_{0}\to0$ and using the continuity of mobile-product the following holds
\begin{multline}
\displaystyle\limsup_{k\to+\infty} \int_{X(\leq s)\cap\Omega_{k,c}}\sum_{j=0}^{s}\binom{n}{j}(-1)^{s-j}w_{k,c}^{n-j}\wedge(T_k+v_{k,c})^{j}\leq\\
\sum_{j=0}^{s}\binom{n}{j}(-1)^{s-j}L^{n-j}\langle F^{j}\rangle
\end{multline}
For the second addendum in (\ref{eq3}), we intersect with $\Omega_{k,c_{1}}$ thus
\begin{multline}\label{mult3}
 \int_{X(\leq s)\cap\Omega_{k,c_{0}}^{c}}(-1)^{s}(\omega-T_k)^{n}=\\
\int_{X(\leq s)\cap\Omega_{k,c_{0}}^{c}\cap\Omega_{k,c_{1}}}(-1)^{s}(\omega-T_k)^{n}+\\
\int_{X(\leq s)\cap\Omega_{k,c_{0}}^{c}\cap\Omega_{k,c_{1}}^{c}}(-1)^{s}(\omega-T_k)^{n}.
\end{multline}
Then for the first addendum in (\ref{mult3}), letting $k\to+\infty$ and $c_{1}\to\nu_{1}^{'}$, we get
\begin{equation*}
 \displaystyle\limsup_{k\to+\infty}\int_{X(\leq s)\cap\Omega_{k,c_{0}}^{c}\cap\Omega_{k,c_{1}}}(-1)^{s}(\omega-T_k)^{n}\leq
\sum_{t=1}^{M_{1}}(L+\nu_{t}^{'}\{u\})^{n-s}\nu_{t}[Y_{t}]
\end{equation*}
 where the sum is over all $Y_{t}\subseteq E_{nn}(F)$ irreducible components of codimension $s$ such that $\nu(T_{min},Y_{t})=\nu_{1}^{'}$.
Now it is clear how to study the second addendum in (\ref{mult3}), one can intersect with $\Omega_{k,c_{2}}$. Then in general we have the following
situation:
\begin{equation}\label{eq4}
 \displaystyle\limsup_{k\to+\infty}\int_{X(\leq s)\cap\left(\bigcap_{i=0}^{j-1}\Omega_{k,c_{i}}^{c}\right)\cap\Omega_{k,c_{j}}}(-1)^{s}(\omega-T_{k})^{n}\leq\sum_{t=M_{j-1}}^{M_{j}}(L+\nu_{t}^{'}\{u\})^{n-s}\nu_{t}[Y_{t}]
\end{equation}
where for all $j=2,\dots,h$, $Y_{t}\subseteq E_{nn}(F)$ are the irreducible components of codimension $s$ with  $\nu(T_{min},Y_t)=\nu_{j}^{'}$
for $M_{j-1}\leq t\leq M_{j}$.\\
While one has that
\begin{equation}
 \int_{X(\leq s)\cap\left(\bigcap_{i=1}^{N}\Omega_{k,c_{i}}^{c}\right)}(-1)^{s}(\omega-T_{k})^{n}\leq\int_{X(\leq s)\cap\Omega_{k,\nu}}(-1)^{s}(\omega-T_{k})^{n}\xrightarrow{k\to+\infty}0
\end{equation}
for dimensional reasons.\\
Then putting togheter the  inequalities of (\ref{eq4}), we finally have
\begin{equation}
 \displaystyle\limsup_{k\to+\infty}\int_{X(\leq s)}(-1)^{s}(\omega-T_k)^{n}\leq\sum_{j=0}^{s}\binom{n}{j}(-1)^{s-j}L^{n-j}\langle F^{j}\rangle+
\sum_{t=1}^{N}\binom{n}{s}(L+\nu_{t}^{'}\{u\})^{n-s}\nu_{t}[Y_t].
\end{equation}
And using now the standard holomorphic Morse inequalities  the statement holds.\\
 If $F$ is not in the interior of $\mathcal{M}^{s}$ and $L$ is not ample,
we consider $F+\varepsilon A$, $L+\varepsilon A$ and then pass to the limit by using the very definition of multicity
(Definition \ref{defini}) .\hfill$\square$
\begin{cor}[Second formulation]
 Let $X$ be a compact projective manifold of complex dimension $n$. Let $L$ and $F$ be two line bundles over $X$ with $L$ nef and
$F\in\mathcal{M}^{s}$ with $\dim (E_{nn}(F))=n-s$. Then we have the following Morse inequalities
\begin{multline}\label{eq2}
 \displaystyle \sum_{j=0}^{s}(-1)^{s-j}h^{j}(X,k(L-F))\leq\\
 \frac{k^{n}}{n!}\sum_{j=0}^{s}\binom{n}{j}(-1)^{s-j}L^{n-j}(F^{j})+\\
 \binom{n}{s}(L+b\{u\})^{n-s}(F^{s}-\langle F^{s}\rangle)+o(k^n)
\end{multline}
where $\displaystyle b=\max_{j=1,\dots, N} \nu_{j}^{'}$.
 \end{cor}
\textit{Proof.} It is sufficient to notice that $b\geq\nu_{j}^{'}$ for all $j=1,\cdots,N$ and 
$F^{s}-\langle F^{s}\rangle=\sum_{j=1}^{N}\nu_{j}[Y_{j}]$, then one can applies Theorem (\ref{Mor1}).\hfill$\square$
\\
\flushleft
Finally, by using remark \ref{rmk9} we can infer that Theorem (\ref{thm5}) is a particular case of Theorem (\ref{Mor1}).
\newpage

\section{Transformation of Lelong numbers by Direct Images}

\subsection{The push-forward of Lelong numbers by a modification}

We want to prove the following:
\begin{prop}\label{pro4}
Let $X$ be a  complex compact manifold with \\ $\dim_{\mathbb{C}}(X)=3$.
Let $\tilde{\mu}:\tilde{X}\rightarrow X$ be a modification of $X$ and $\Omega$ is a smooth, positive form on $\tilde{X}$ of bidimension 1. Then
$\nu(\tilde{\mu_{*}}(\tilde{\Omega}),Y)=0$ $\forall Y$  irreducible curve on $X$, where $\nu$ is the generic Lelong number.
\end{prop}

Since a modification can be decomposed as a finite sequence of blow-ups with smooth centers, we can write $\tilde{\mu}$ as follows:
\begin{equation}
 \tilde{\mu}:\tilde{X}:=X_{s}\rightarrow X_{s-1}\rightarrow\cdots\rightarrow X_{1}\rightarrow X_{0}=X. 
\end{equation}
with $\forall i=1,\dots,s\quad Z_{i},E_{i}\subset X_{i}$ are the smooth centers and the exceptional divisors of the $i$-th blow-up respectively.
\\
Let $y\in Y$ be a generic point, without loose of generality one can suppose that $y$ is the origin for a local chart with coordinates
$z=(z_{1},z_{2},z_{3})$. One has the following:
\begin{eqnarray*}
 \nu(\tilde{\mu_{*}}(\tilde{\Omega}),0) & = &\lim_{r\rightarrow0}\dfrac{1}{r^{2}}\int_{B(0,r)}\tilde{\mu_{*}}(\tilde{\Omega})(z)\wedge\dfrac{i}{2\pi}\partial\bar{\partial}|z|^{2}= \nonumber \\
&& = {}\lim_{r\rightarrow0}\dfrac{1}{r^{2}}\int_{\tilde{\mu}^{-1}(B(0,r))}\tilde{\Omega}\wedge\tilde{\mu^{*}}\left(\dfrac{i}{2\pi}\partial\bar{\partial}|z|^{2}\right).
\end{eqnarray*}
Let's note that the integrand in the last equality is smooth so one wants to estimate the volume of $\tilde{\mu}^{-1}(B(0,r))$ and to compare with the parameter $r^{2}$.\\
Using the fact that $\tilde{\mu}$ can be decomposed as a finite sequence of blow-ups, $\tilde{\mu}^{-1}(B(0,r))$ can be expressed explicitly by a local 
expression of the composition of these blow-ups. Thus, one can associate to this local expression a $psh$-function $\varphi$ such that the volume of
$\tilde{\mu}^{-1}(B(0,r))$ can be studied in therms of asymptotic estimates for the volume of sublevel sets $\{\varphi<\text{log}\;r\}$.\\
Explicitly, let $B(0,r)=\{|z_{1}|^{2}+|z_{2}|^{2}+|z_{3}|^{2}<r\}$ be the local expression on $X$ with coordinates $(z_{1},z_{2},z_{3})$, so
$\tilde{\mu}^{-1}(B(0,r))=\{|h_{1}|^{2}+|h_{2}|^{2}+|h_{3}|^{2}<r\}$ where $h_{i}:=h_{i}(\tilde{z_{1}},\tilde{z_{2}},\tilde{z_{3}})$ are holomorphic
functions on $\tilde{X}$ and $(\tilde{z_{1}},\tilde{z_{2}},\tilde{z_{3}})$ are local coordinates on $\tilde{X}$.\\
Then one defines $\displaystyle \varphi=\frac{1}{2}\text{log}\sum_{i=1}^{3}|h_{i}|^{2}$, so that $\tilde{\mu}^{-1}(B(0,r))=\{\varphi<\text{log}\;r\}.$ \\
Let us give the following 
\begin{defn}[\cite{DK}]\label{def1}
 Let $K\subset X$ be a compact set, $U\Subset X$ a relatively compact neighborhood of $K$ and let $\theta_{U}$ be the Lebesgue measure on $U$ 
associated with some choice of hermitian metric $\omega$ on $X$.\\ Then the \textbf{log-canonical threshold} of $\varphi$ is defined as:
\begin{equation*}
 c_{K}(\varphi):=\sup~\{c\geq0~:r^{-2c}\theta_{U}(\{\varphi<\text{log}~r\})~\text{is bounded as}~r\to0~\text{for}~U\supset K\}.
\end{equation*}
\end{defn}
Now the crucial fact is that each $h_{i}$ can be expressed as a sum of a holomorphic monomial in the $\tilde{z_{i}}$ and another holomorphic
function, so that looking at $\tilde{\mu}(\tilde{z_{1}},\tilde{z_{2}},\tilde{z_{3}})=(h_{1}(\tilde{z_{1}},\tilde{z_{2}},\tilde{z_{3}}),h_{2}(\tilde{z_{1}},\tilde{z_{2}},\tilde{z_{3}}),h_{3}(\tilde{z_{1}},\tilde{z_{2}},\tilde{z_{3}}))$
as an ideal ``almost'' monomial, one can calculate $c_{K}(\varphi)$ using Howald's theorem \cite{How1}.
\\
\begin{rmk}
 One may assume that for all $i$ the centers $Z_{i}$ are connected by considering a longer sequence of blow-ups if necessary.
\end{rmk}
\begin{rmk}\label{rmk3}
 If $Y$ is not entirely contained in $Z_{0}$ then the statement is trivial. In fact let's suppose $Y\cap Z_{0}=\{a_{1},\dots,a_{k}\}$ is a finite
set of points or possibly  the empty set, let $\mu_{1}\circ\dots\circ\mu_{s}:X_{s}\to X_{s-1}\to\cdots X_{1}\to X_{0}$ be a finite sequence
of blow-ups, let $\tilde{\Omega}$ be a smooth positive form of bidimension 1 on $X_{s}$, and let $E_{1}$ be the exceptional divisor of $\mu_{1}$
then $\mu_{1}(X_{1}\setminus E_{1})=X_{0}\setminus\mu_{1}(E_{1})=X_{0}\setminus Z_{0}$ and $\mu_{1}|_{X_{1}\setminus E_{1}}:X_{1}\setminus E_{1}\to X_{0}\setminus Z_{0}$ is
an isomorphism. Thus $\nu(\tilde{\mu}_{*}(\tilde{\Omega}),Y)=
\nu((\mu_{2}\circ\dots\circ\mu_{s})_{*}\tilde{\Omega},\mu_{1}^{-1}(Y))$, and now if $\mu_{1}^{-1}(Y)$ is not entirely contained in $Z_{1}$, the center of $\mu_{2}$,
 one can omit $\mu_{2}$ in the calculation of $\nu$ and so on. So we may assume $Y=Z_{0}$. 
\end{rmk}
Let $\tilde{\mu}=\mu_{1}\circ\dots\circ\mu_{s}$ be a finite sequence of blow-ups with smooth connected centers $Z_{i}$ for $=1,\dots,s-1$
( $Z_{0}$ is given by (\ref{rmk3})). Let's define $Exc(\mu_{1}\circ\dots\circ\mu_{i}):=\cup_{j=1}^{i}(\mu_{j}\circ\dots\circ\mu_{i})^{-1}(Z_{j-1})$ 
 then the centers can just satisfy one of these conditions:
\item \textbf{a)} $Z_{i}$ is a point and $Z_{i}\in Exc(\mu_{1}\circ\dots\circ\mu_{i})$,\\
\item \textbf{b)} $Z_{i}$ is a point and $Z_{i}\notin Exc(\mu_{1}\circ\dots\circ\mu_{i})$,\\
\item \textbf{c)} $Z_{i}$ is a curve such that $Z_{i}\cap Exc(\mu_{1}\circ\dots\circ\mu_{i})=\emptyset$,\\
\item \textbf{d)} $Z_{i}$ is a curve such that $Z_{i}\cap Exc(\mu_{1}\circ\dots\circ\mu_{i})$ is a finite set of points,\\
\item \textbf{e)} $Z_{i}$ is a curve such that $Z_{i}\subset Exc(\mu_{1}\circ\dots\circ\mu_{i})$ entirely.\\

\newpage Then one can associate to this sequence a rooted tree $\mathcal{T}$, where the root is 
the vertex $X_{0}$ and the edges $\mu_{i}^{a},\mu_{i}^{b},\mu_{i}^{c},\mu_{i}^{d},\mu_{i}^{e}$ denote the 
$i$-th blow-up (for $i\geq2$) where the centers $Z_{i}$ satisfy respectively \textbf{a),b),c),d),e)}.\\
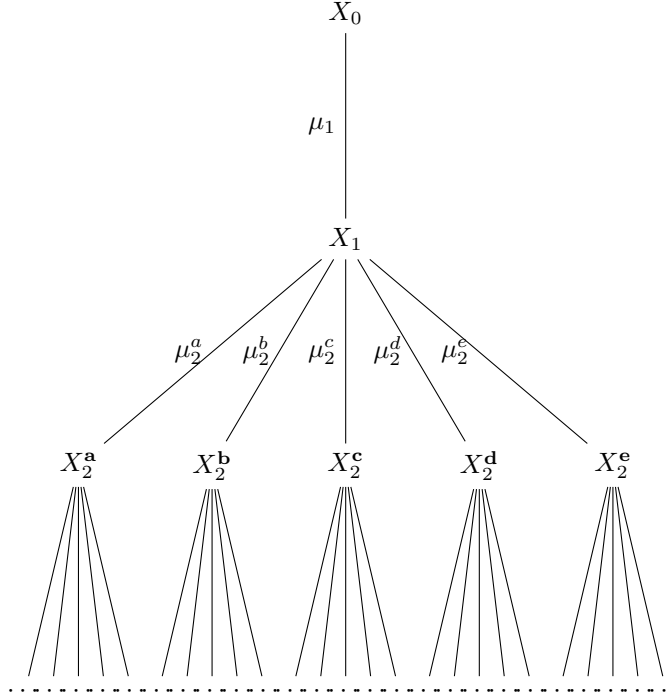
\begin{center}
 \begin{tikzpicture}
[ level 1/.style={sibling distance=6em},
level 2/.style={sibling distance=5em}, level distance=1cm,
level 3/.style={sibling distance=1em}, level distance=3cm
]
  \node {$X_{0}$}
child { node {$X_{1}$}
child{ node {$X_{2}^{\textbf{a}}$}
child { node{$\cdots$}} child{ node{$\cdots$}} child{ node{$\cdots$}} child{ node{$\cdots$}} child{ node{$\cdots$}} edge from parent node[left] {$\mu_{2}^{a}$} }
child{ node {$X_{2}^{\textbf{b}}$}
child { node{$\cdots$}} child{ node{$\cdots$}} child{ node{$\cdots$}} child{ node{$\cdots$}} child{ node{$\cdots$}} edge from parent node[left] {$\mu_{2}^{b}$} }
child{ node {$X_{2}^{\textbf{c}}$}
child { node{$\cdots$}} child{ node{$\cdots$}} child{ node{$\cdots$}} child{ node{$\cdots$}} child{ node{$\cdots$}} edge from parent node[left] {$\mu_{2}^{c}$} }
child{ node {$X_{2}^{\textbf{d}}$}
child { node{$\cdots$}} child{ node{$\cdots$}} child{ node{$\cdots$}} child{ node{$\cdots$}} child{ node{$\cdots$}} edge from parent node[left] {$\mu_{2}^{d}$} }
child{ node {$X_{2}^{\textbf{e}}$}
child { node{$\cdots$}} child{ node{$\cdots$}} child{ node{$\cdots$}} child{ node{$\cdots$}} child{ node{$\cdots$}} edge from parent node[left] {$\mu_{2}^{e}$} }
edge from parent 
node[left] {$\mu_{1}$}};
 \end{tikzpicture}
\end{center}
Since we are interested in calculating $\nu(\tilde{\mu_{*}}\tilde{\Omega},Z_{0}),$ where $\tilde{\mu}=\mu_{1}\circ\dots\circ\mu_{s}$ is a finite
sequence of blow-ups with smooth centers, we want to understand which path of the tree $\mathcal{T}$ represents the more general sequence according
to the following
\begin{defn}\label{def2}
The sequence of blow-ups $\mu_{1}\circ\dots\circ\mu_{s}$ is called minimal if there does not exist another sequence $\mu_{1}^{'}\circ\dots\circ\mu_{s}^{'}$
such that $\nu((\mu_{1}\circ\dots\circ\mu_{s})_{*}(\tilde{\Omega}),Z_0)=\nu((\mu_{1}^{'}\circ\dots\circ\mu_{s}^{'})_{*}(\tilde{\Omega}),Z_0)$ where 
there exists an $i$ for which $\mu_{i}^{'}$ is a local biholomorphism.
\end{defn}
Thus we want to understand which path of the tree doesn't have this property because it reflects the general sequence of $s$ blow-ups. The following
lemma shows  that if $Z_{i}$ satisfies one of conditions \textbf{a),b),c),d)} then the corresponding path is not minimal.
\begin{lem}\label{lem1}
 Let $\mathcal{T}$ be the rooted tree associated to a finite sequence of blow-ups. The following path
\begin{equation*}
 \tilde{\mu}^{J_{1},J_{2},\dots,J_{s}}:=\mu_{1}^{J_{1}}\circ\mu_{2}^{J_{2}}\circ\cdots\circ\mu_{s}^{J_{s}}~\textit{where}~J_{1},\dots,J_{s}\in\{\textbf{a,b,c,d}\}
\end{equation*}
is not minimal.
\end{lem}
\textit{Proof.} 
\textbf{Second blow-up.}  Let $\tilde{\mu}=\mu_{1}^{J}\circ\mu_{2}\circ\mu_{3}\circ\cdots\circ\mu_{s}$ be a sequence of $s$ blow-ups such that $J\in\{a,b,c,d\}.$\\
\emph{\underline{Case J=\textbf{a}.}} Let $p_{0}\in Z_{0}$ such that $\mu_{1}^{-1}(p_0)\ni Z_1$ then for all $p\neq p_{0}$ with $p\in Z_{0}$ there exists $r>0$
small enough with $Z_{1}\nsubseteq\mu_{1}^{-1}(B(p,r))$. Then near $p$
\begin{equation}\label{eq5}
 \mu_{2}:\mu_{2}^{-1}(\mu_{1}^{-1}(B(p,r)))\xrightarrow{\simeq}\mu_{1}^{-1}(B(p,r))
\end{equation}
is an isomorphism. Now if we consider $\tilde{\mu}^{'}=\mu_{1}\circ\mu_{2}^{'}\circ\mu_{3}\circ\cdots\circ\mu_{s}$ with $\mu_{2}^{'}$ any isomorphism, as above we have that
$\nu(\tilde{\mu_{*}}\tilde{\Omega},Z_{0})=\nu(\tilde{\mu}_{*}^{'}\tilde{\Omega},Z_{0})$.\\
\emph{\underline{Case J=\textbf{b}.}}  The  argument above  is still valid.\\
\emph{\underline{Case J=\textbf{c}}} For all points $p\in Z_{0}$ one can find an $r>0$ sufficiently small such that (\ref{eq1}) still holds,
then one can repeat the argument of case $J=\textbf{a}$.\\
\emph{\underline{ Case J=\textbf{d}}}. Let $\{a_{1},\dots,a_{k}\}=Z_{1}\cap E_{1}=Exc(\mu_{1})$ be a finite set of points.
 For all $p\in Z_{0}$ such that $p\neq\mu_{1}(a_{i})$ for all $i=1,\dots k$ one can find
again an $r>0$ small enough such that $a_{i}\notin\mu_{1}^{-1}(B(p,r))$ $\forall i=1,\dots,k$. Repeating the argument as before one can find another
sequence of $s-1$ blow-ups $\tilde{\mu}^{'}$ such that $\nu(\tilde{\mu}_{*}\tilde{\Omega},Z_{0})=\nu(\tilde{\mu}_{*}^{'}\tilde{\Omega},Z_{0})$.\\
\textbf{Third blow-up.} Now let's suppose $\tilde{\mu}=\mu_{1}^{e}\mu_{2}^{J}\circ\mu_{3}\circ\cdots\circ\mu_{s}$ with $J\in\{\textbf{a,b,c,d}\}$ and we may assume that $\mu_{1}:Z_{1}\to Z_{0}$ is surjective.\\
\emph{\underline{Case J=\textbf{a,b}}}. Let $p$ be a point of $Z_{0}$ such that $(\mu_{1}\circ\mu_{2})(Z_{2})=p$ then for all $p^{'}\neq p$ and $p^{'}\in Z_{0}$
we can find a small enough $r>0$ such that
\begin{equation}\label{eq6}
\mu_{3}:\mu_{3}^{-1}\left(\mu_{2}^{-1}\circ\mu_{1}^{-1}(B(p^{'},r))\right)\xrightarrow{\simeq}\mu_{2}^{-1}\circ\mu_{1}^{-1}(B(p^{'},r)).
\end{equation}
Then for any isomorphism $\mu_{3}$ as above one can change the sequence $\tilde{\mu}$ with  $\tilde{\mu}^{'}=\mu_{1}^{e}\circ\mu_{2}^{J}\circ\mu_{3}^{'}\circ\cdots\circ\mu_{s}$
and we have $\nu(\tilde{\mu}_{*}\tilde{\Omega},Z_{0})=\nu(\tilde{\mu}_{*}^{'}\tilde{\Omega},Z_{0})$.\\
\emph{\underline{Case J=\textbf{c}}}. The argument above is still valid.\\
\emph{\underline{Case J=\textbf{d}}}. Let $\{a_{1},\cdots,a_{k}\}=Z_{2}\cap Exc(\mu_{1}^{e}\circ\mu_{2})$ be a finite set of points in $X_{2}^{\textbf{d}}$, then
there exist $p_{1},\dots,p_{k}\in Z_{0}$ such that $(\mu_{1}^{e}\circ\mu_{2}^{d})(a_{i})=p_{i}$. 
Then for all $p^{'}\in Z_{0}$ with $p^{'}\neq p_{i}$ $1\leq i\leq k$ we can find $r>0$ sufficiently small so that (\ref{eq6}) is again true.
Therefore one finds a sequence of $s-1$ blow-ups as in the cases above.
\\
Now it's clear how to study the general case.\\
\textbf{i-th blow-up.} Let $\tilde{\mu}=\mu_{1}^{e}\circ\dots\circ\mu_{i-1}^{e}\circ\mu_{i}^{J}\circ\mu_{i+1}\circ\cdots\circ\mu_{s}$ be a finite sequence
of blow-ups with $J\in\{\textbf{a,b,c,d}\}.$\\
\emph{\underline{Case J=\textbf{a}}}. Let $Z_{i}$ be the center of $\mu_{i+1}$ with $Z_{i}\subset Exc(\mu_{1}^{e}\circ\cdots\circ\mu_{i-i}^{e}\circ\mu_{i}^{a}).$
and let 
$(\mu_{1}^{e}\circ\cdots\circ\mu_{i-i}^{e}\circ\mu_{i}^{a})(Z_{i})=p$. Thus for all $p^{'}\neq p$ in $Z_{0}$ we can find a small positive $r$ such that the following hold
\begin{itemize}
 \item  1) $(\mu_{1}^{e}\circ\cdots\circ\mu_{i}^{a})^{-1}(B(p^{'},r))\cap Z_{i}=\emptyset$ \\
\item   2) $\mu_{i+1}:\mu_{i+1}^{-1}\left((\mu_{1}^{e}\circ\cdots\circ\mu_{i-i}^{e}\circ\mu_{i}^{a})^{-1}B(p^{'},r)\right)\xrightarrow{\simeq}(\mu_{1}^{e}\circ\cdots\circ\mu_{i-i}^{e}\circ\mu_{i}^{a})^{-1}B(p^{'},r)$.
\end{itemize}
So $\tilde{\mu}^{'}=\mu_{1}^{e}\circ\cdots\mu_{i-1}^{e}\circ\mu_{i}^{a}\circ\mu_{i+1}^{'}\circ\cdots\mu_{s}$, where $\mu_{i+1}^{'}$ is an isomorphism,
satisfies $\nu(\tilde{\mu}_{*}\tilde{\Omega},Z_{0})=\nu(\tilde{\mu}_{*}^{'}\tilde{\Omega},Z_{0})$ .\\
\emph{\underline{Case J=\textbf{b}}}. As in the case J=\textbf{a} 1) and 2) are still valid,but now they are valid for all $p^{'}\in Z_{0}$.\\
\emph{\underline{Case J=\textbf{c}}}. Again as in the case J=\textbf{a,b} we can find an $r>0$ sufficiently small such that for all $p^{'}\in Z_{0}$
1) and 2) are still valid and then the statement holds.\\
\emph{\underline{Case J=\textbf{d}}}.  $Z_{i}\cap Exc(\mu_{1}^{e}\circ\cdots\circ\mu_{i-i}^{e}\circ\mu_{i}^{a})=\{a_{1},\dots,a_{k}\}$,
where $a_{i}$ are points of $Z_{i}.$ Now there exist $p_{1},\dots,p_{k}\in Z_{0}$ such that $(\mu_{1}^{e}\circ\cdots\circ\mu_{i-i}^{e}\circ\mu_{i}^{a})(a_{i})p_{i}$.
So we can find $r>0$ small so that  for $p\neq p_{i}$ $1\leq i\leq k$ such that
\begin{equation*}
a_{i}\notin (\mu_{1}^{e}\circ\cdots\circ\mu_{i-i}^{e}\circ\mu_{i}^{a})^{-1}B(p,r),
\end{equation*}
so we can repeat the same reasoning as in the previous cases.\hfill$\square$
\\

According to Lemma (\ref{lem1}) it remains to analyze just the path of the tree $\mathcal{T}$ corresponding to a sequence of blow-ups
 $\tilde{\mu}=\mu_{1}^{e}\circ\cdots\circ\mu_{s-1}^{e}\circ\mu_{s}^{e}$ where for all $i=1,\dots,s-1$ $Z_{i}\subset Exc(\mu_{1}^{e}\circ\cdots\circ\mu_{i}^{e})$ and
$\mu_{1}^{e}:Z_{1}\to Z_{0}$ is  surjective. Now let's distinguish two possible cases:
\begin{equation}\label{eq7}
 1)~ Z_{i}\subseteq\widetilde{Exc}(\mu_{1}^{e}\circ\cdots\circ\mu_{i}^{e}):=Exc(\mu_{1}^{e}\circ\cdots\circ\mu_{i}^{e})\setminus E_{i};\quad 2)~ Z_{i}\subseteq E_{i}~\textit{entirely.}
\end{equation}

\begin{nota}
 We write
\begin{equation}
 \tilde{\mu}^{e}=\tilde{\mu}_{1}^{e_{i_{1}}}\circ\tilde{\mu}_{2}^{e_{i_{2}}}\cdots\circ\tilde{\mu}_{s}^{e_{i_{s}}}
\end{equation}
with  $i_{1},\dots,i_{s}\in\{1,2\}$, where $1,2$ correspond to the two conditions of (\ref{eq7}) for the center of the following blow-up.
(For example $\mu_{j}^{e_{1}}$ means that the center of the $j+1$-th blow-up satisfies condition 1 of (\ref{eq7}).)
\end{nota}
Let's note that for the first blow-up it is not necessary to indicate which condition of (\ref{eq7}) $Z_{1}$ satisfies thanks to Lemma \ref{lem1} and to the following remark.\\
\begin{rmk}\label{rmk5}
 If $Z_{i}\subseteq E_{i}$ for some $i=1,\dots,s-1.$ Since
\begin{equation*}
 \mu_{i}|_{E_{i}}:E_{i}\to Z_{i-1}
\end{equation*}
is a holomorphic fiber bundle isomorphic to the projectivized normal bundle $\mathbb{P}(N_{Z_{i-1}})\to Z_{i-1}.$ We can suppose that
$\mu_{i}|_{Z_{i}}:Z_{i}\to Z_{i-1}$ is surjective.
\end{rmk}

\flushleft Now let's fix some notations. Since everything is local, without loose of generality, one can suppose that $Z_{0}=\{x_{0}=y_{0}=0\}$ where $(x_{0},y_{0},z_{0})$
are local coordinates on $X_{0}=X$.
\begin{rmk}\label{rmk6}
 Let Y be a $n$-dimensional smooth variety and $B$ a closed smooth subvariety with codim$_{Y}(B)=t$ and $\sigma:\tilde{Y}\rightarrow Y$ the blow-up
of $Y$ with center $B$. \\ For every point $b_{0}\in B$ there exists local coordinates $(u_{1},\dots,u_{n})$ on $Y$ centered at $b_{0}$ such that
$B=\{u_{1}=\dots=u_{t}=0\}$ and local coordinates $(w_{1},\dots,w_{n})$ on $\tilde{Y}$ such that the map $\sigma$ is given by:
\begin{equation}
 \sigma(w_{1},\dots,w_{n})=(w_{1}w_{j},\dots,w_{j-1}w_{j};w_{j};w_{j+1}w_{j},\dots,w_{t}w_{j},w_{t+1},\dots,w_{n})
\end{equation}
$\forall j=1,\dots,t.$
\end{rmk}
\newpage
\begin{lem}
 If the first condition of (\ref{eq7}) is satisfied for some $i$ with $1\leq i\leq s$ the $\tilde{\mu}$ is not minimal.
\end{lem}
\textit{Proof.}
Let's consider the $i$-th blow-up $\mu_{i}=X_{i}\to X_{i-1}$ by using (\ref{rmk6}) we see that there exists local coordinates $(x_{i},y_{i},z_{i})$ on $X_{i}$
such that
\begin{equation*}
 \mu_{i}(x_{i},y_{i},z_{i})=\left\{\begin{array}{ll}
                                    (x_{i},y_{i}x_{i},z_{i})\\
                                    (x_{i}y_{i},y_{i},z_{i})
                                   \end{array}\right.
\end{equation*}
where we used the same coordinates on two different charts to keep notations simple. By Lemma (\ref{lem1}) one has that $Z_{i}\subseteq Exc(\mu_{1}^{e}\circ\cdots\circ\mu_{i}^{e})$.
Now let us distinguish the following two cases:
\begin{enumerate}
 \item $Z_{i}\subseteq E_{i}$ entirely,\\
\item $Z_{i}\subseteq\widetilde{Exc}(\mu_{1}^{e}\circ\cdots\circ\mu_{i}^{e})$ entirely.
\end{enumerate}
\underline{\emph{Case 1.}} Repeating the argument of (\ref{rmk5}) we get that
\begin{equation*}
 Z_{i}=\{x_{i}=0=f_{i}(y_{i},z_{i})\}~\text{or}~Z_{i}=\{y_{i}=0=g_{i}(x_{i},z_{i})\}
\end{equation*}
depending on which chart we are considering where $f_{i},g_{i}$ are holomorphic functions. By surjectivity of the map $\mu_{i}:Z_{i}\to Z_{i-1}$
we deduce that on the generic point of $Z_{i}$
\begin{equation*}
 \frac{\partial f_i}{\partial y_i}\neq0\quad\text{or}\quad\frac{\partial g_i}{\partial x_i}\neq0
\end{equation*}
then by implicit function theorem
there exist holomorphic functions $F_{i},G_{i}$ such that
\begin{equation*}
 Z_{i}=\{x_{i}=0,y_{i}=F{z_{i})}\}~\text{or}~Z_{i}=\{y_{i}=0,x_{i}=G_{i}(z_{i})\}.
\end{equation*}
Now considering the following holomorphic changes of coordinates on the two different charts of $\mu_{i}$
\begin{equation}\label{eq8}
 \left\{\begin{array}{l}
         U_{i}=x_{i}\\
         V_{i}=y_{i}-F_{i}(z_{i})\\
         W_{i}=z_{i}
        \end{array}\right.
\quad\text{or}~\left\{\begin{array}{l}
             U_{i}=x_{i}-G_{i}(z_{i})\\
             V_{i}=y_{i}\\
             W_{i}=z_{i}
            \end{array}\right.
\end{equation}
where we have still used the same coordinates on the two different charts. So the equation of $Z_{i}$ by (\ref{eq8}) becomes $Z_{i}=\{U_{i}=V_{i}=0\}$.
 Now we can use (\ref{rmk6}) blowing-up along $Z_{i}$ and repeating the same operations as above.\\

\flushleft\underline{\emph{Case 2.}} Let $i\in\{1,\dots,s-1\}$ be the first index such that $Z_{i}\nsubseteq E_{i}$ entirely, i.e. $Z_{i}\cap E_{i}$ is a finite
set of points. Hence for all $j<i$ $Z_{j}\subseteq E_{j}$ entirely. Now to keep notations simple we denote by $E_{j}^{'}$ for $j<i$ the pre-image
of $E_{j}\setminus Z_{j}$ in $X_{i}$ by $(\mu_{j}\circ\cdots\circ\mu_{i})$; then $\displaystyle Z_{i}\subseteq\bigcup_{j=1}^{i-1}E_{j}^{'}$.
Since $Z_{i}$ is connected then we can suppose that $Z_{i}\subseteq E_{j_{0}}$ entirely for some $1\leq j_{0}\leq i-1$.
The expression in local coordinates of $\mu_{i}$ is given by
\begin{equation*}
 \mu_{i}(x_i,y_i,z_i)=\left\{\begin{array}{l}
                              (x_{i},y_{i}x_{i},z_{i})\\
                              (x_{i}y_{i},y_{i},z_{i})
                             \end{array}\right.
\end{equation*}
where $(x_i,y_i,z_i)$ depends on $(U_{i-1},V_{i-1},W_{i-1})$ on $X_{i-1}.$ But now $Z_{i}\nsubseteq E_{i}$ entirely means that $(\mu_{j_{0}+1}\circ\cdots\circ\mu_{i})|_{Z_{i}\setminus\{p_{1},\dots,p_{k}\}}$ is 
an isomorphism and furthermore $(\mu_{j_{0}+1}\circ\cdots\circ\mu_{i})(Z_{i}\setminus\{p_{1},\dots,p_{k}\})\subseteq E_{j_{0}}$ entirely thus in local coordinates $(x_{j_{0}},y_{j_{0}},z_{j_{0}})$
on $X_{j_{0}}$ we have that
\begin{equation*}
 \mu_{j_{0}}(x_{j_{0}},y_{j_{0}},z_{j_{0}})=\left\{\begin{array}{l}
                              (x_{j_{0}},y_{j_{0}}x_{j_{0}},z_{j_{0}})\\
                              (x_{j_{0}}y_{j_{0}},y_{j_{0}},z_{j_{0}})
                             \end{array}\right.
\end{equation*}
and $(\mu_{j_{0}+1}\circ\cdots\circ\mu_{i})(Z_{i})$ is equal to
 \begin{equation}\label{eq9}
  \{x_{j_{0}}=0=f_{i}(y_{j_{0}},z_{j_{0}})\}~\text{or}~\{y_{j_{0}}=0=g_{i}(x_{j_{0}},z_{j_{0}})\}.
 \end{equation}
So performing the $(i+1)$-blow-up along $Z_{i}$ on $X_{i}$ or along $(\mu_{j_{0}+1}\circ\cdots\circ\mu_{i})(Z_{i})$ on $X_{j_{0}}$ it is essentially equivalent.
By remark (\ref{rmk5}) it is not restrictive to suppose that
\begin{equation*}
 \mu_{j_{0}}|_{(\mu_{j_{0}+1}\circ\cdots\circ\mu_{i})(Z_{i})}:(\mu_{j_{0}+1}\circ\cdots\circ\mu_{i})(Z_{i})\twoheadrightarrow Z_{j_{0}-1}
\end{equation*}
and this means that in the local expressions (\ref{eq9}) one can suppose that 
\begin{equation*}
 \frac{\partial f_{i}}{\partial y_{j_{0}}}\neq0~\text{or}~\frac{\partial g_{i}}{\partial x_{j_{0}}}\neq0
\end{equation*}
depending which local chart we are considering.\\
So we can define new coordinates on $X_{j_{0}}$ using (\ref{eq8}). Namely:
\begin{equation}
 \left\{\begin{array}{l}
         U_{i}=x_{j_{0}}\\
         V_{i}=y_{j_{0}}-F_{i}(z_{j_{0}})\\
         W_{i}=z_{j_{0}}
        \end{array}\right.
\quad\text{or}\quad
\left\{\begin{array}{l}
        U_{i}=x_{j_{0}}-G_{i}(z_{_{j_{0}}})\\
        V_{i}=y_{j_{0}}\\
        W_{i}=z_{j_{0}}
       \end{array}\right.
\end{equation}
 for which $(\mu_{j_{0}+1}\circ\cdots\circ\mu_{i})(Z_{i})=\{U_{i}=V_{i}=0\}$ 
and we can now perform the $(i+1)$-blow-up. The relation with the coordinates system $(U_{j_{0}},V_{j_{0}},W_{j_{0}})$ on $X_{j_{0}}$ is straightforward by using:
\begin{equation*}
 \left\{\begin{array}{l}
         U_{j_{0}}=x_{j_{0}}\\
         V_{j_{0}}=y_{j_{0}}-F_{j_{0}}(z_{j_{0}})\\
         W_{j_{0}}=z_{j_{0}}
        \end{array}\right.
\quad\text{or}\quad
\left\{\begin{array}{l}
        U_{j_{0}}=x_{j_{0}}-G_{j_{0}}(z_{j_{0}})\\
        V_{j_{0}}=y_{j_{0}}\\
        W_{j_{0}}=z_{j_{0}}.
       \end{array}\right.
\end{equation*}
Then the expression in local coordinates of the composition $(\mu_{j_{0}}\circ\cdots\circ\mu_{i+1})$ is the following:
\begin{equation}
 (\mu_{j_{0}}\circ\cdots\circ\mu_{i+1})=\left\{\begin{array}{l}
                                                (x_{i+1},x_{i+1}(x_{i+1}y_{i+1}F_{i}(z_{i+1})),z_{i+1})\\
                                                (x_{i+1}y_{i+1},x_{i+1}y_{i+1}(y_{i+1}+F_{i}(z_{i+1})),z_{i+1})\\
                                                ((x_{i+1}+G_{i}(z_{i+1}))x_{i+1}y_{i+1},x_{i+1}y_{i+1},z_{i+1})\\
                                                ((x_{i+1}y_{i+1}+G_{i}(z_{i+1}))y_{i+1},y_{i+1},z_{i+1}) 
                                               \end{array}\right.
\end{equation}
Let's observe that in the \underline{\emph{Case 2}} we can disreguard the intersection of $Z_{i}$ with $E_{i}$ since it has
dimension zero and it has no effect on the generic Lelong number. \hfill$\square$                                                                                                                            
                                             
\subsection{Estimate for $c_{K}(\varphi)$}

Let's now compute the log-canonical threshold of the ideals rising from $\tilde{\mu}=\mu_{1}\circ\cdots\circ\mu_{s}$ a composition of $s$ blow-ups.
By the  previous section without loss of generality one can suppose that each center $Z_{i}$ is entirely contained in the exceptional divisor $E_{i}$, since 
 we have already notice (\underline{\emph{ Case 2}}) that the expression of the composition of  blow-ups is equivalent to
a composition of blow-ups in which each center $Z_{i}$ is entirely contained in $E_{i}.$
\begin{prop}\label{pro5}
 Let $\tilde{\mu}=\mu_{1}\circ\cdots\circ\mu{s}:\tilde{X}\rightarrow X$ be a finite sequence of blow-ups with smooth centers $Z_{i}$ such that
$\dim_{\mathbb{C}}(Z_{i})=1$ and $Z_{i}$ is entirely contained in the exceptional divisor of $\mu_{i}$ and it's surjective $\forall i=1,\dots, s$.
Then considering  the ideal generated by all possible $2^{k}$ expressions of $\tilde{\mu}$ in local coordinates, one can reduce them to ideals of the forms
\begin{equation*}
\mathcal{I}=(x_{s}^{h}y_{s}^{k},z_{s})~\text{with}~ h,k\geq0~\text{with}~max(h,k)>0.
\end{equation*}
\end{prop}
\textit{Proof.} Finite induction on $s$.\\
\textbf{Case s=1.} $\tilde{\mu}=\mu_{1}:X_{1}\rightarrow X_{0}$, there exist local coordinates $(x_{0},y_{0},z_{0})$ on $X_{0}$ such that $Z_{0}=\{x_{0}=y_{0}=0\}$, by (\ref{rmk6}) the expressions of $\mu_{1}$ in local coordinates are:
\begin{equation*}
 \mu_{1}(x_{1},y_{1},z_{1})=(x_{1},x_{1}y_{1},z_{1})~\text{and}~\mu_{1}(x_{1},y_{1},z_{1})=(x_{1}y_{1},y_{1},z_{1})
\end{equation*}
so one has that
\begin{equation*}
 (x_{1},x_{1}y_{1},z_{1})=(x_{1},z_{1})~\text{and}~(x_{1}y_{1},y_{1},z_{1})=(y_{1},z_{1}).
\end{equation*}
\textbf{Case s=2.} $\tilde{\mu}=\mu_{1}\circ\mu_{2}:X_{2}\xrightarrow{\mu_{2}}X_{1}\xrightarrow{\mu_{1}}X_{0}.$\\
Let $Z_{0}=\{x_{0}=y_{0}=0\}$ be the center of $\mu_{1}$ then there exists local coordinates on $X_{1}$ such that
\begin{equation*}
 \mu_{1}(x_{1},y_{1},z_{1})=(x_{1},x_{1}y_{1},z_{1})~\text{and}~\mu_{1}(x_{1},y_{1},z_{1})=(x_{1}y_{1},y_{1},z_{1})
\end{equation*}
Let's consider first the expression $\mu_{1}(x_{1},y_{1},z_{1})=(x_{1},x_{1}y_{1},z_{1})$; by (\ref{eq8}) let $Z_{1}=\{x_{1}=0=f_{1}(y_{1},z_{1})\}=\{x_{1}=0,~y_{1}=F_{1}(z_{1})\}$
be the center of $\mu_{2}$.  So that $Z_{1}=\{U_{1}=V_{1}=0\}.$
Then there exist local coordinates on $X_{2}$ such that:
\begin{equation*}
  \mu_{2}(x_{2},y_{2},z_{2})=(x_{2},x_{2}y_{2},z_{2})~\text{and}~\mu_{2}(x_{2},y_{2},z_{2})=(x_{2}y_{2},y_{2},z_{2}),
\end{equation*}
repeating the argument used in the case of $\mu_{1}(x_{1},y_{1},z_{1})=(x_{1}y_{1},y_{1},z_{1})$, the following holds:
\begin{enumerate}
 \item $(\mu_{1}\circ\mu_{2})(x_{2},y_{2},z_{2})=(x_{2},x_{2}(x_{2}y_{2}+F_{1}(z_{2})),z_{2})$\\
\item $(\mu_{1}\circ\mu_{2})(x_{2},y_{2},z_{2})=(x_{2}y_{2},x_{2}y_{2}(y_{2}+F_{1}(z_{2})),z_{2})$\\
\item $(\mu_{1}\circ\mu_{2})(x_{2},y_{2},z_{2})=((x_{2}+G_{1}(z_{2}))x_{2}y_{2},x_{2}y_{2},z_{2})$\\
\item $(\mu_{1}\circ\mu_{2})(x_{2},y_{2},z_{2})=((x_{2}y_{2}+G_{1}(z_{2}))y_{2},y_{2},z_{2}).$
\end{enumerate}
Considering the simmetry of these expressions, it's enough to analyze  1) and 2), thus these local expressions can be reduced to the ideals respectively:
\begin{equation*}
 (x_{2},z_{2})~\text{and}~(x_{2}y_{2},z_{2}).
\end{equation*}
\textbf{Case s=3.} One would have to consider the $2^{3}$ possible expressions in local coordinates but by simmetry it's enough to study  four of them, in particular:
\item $\mathfrak{a}=(x_{3},x_{3}((x_{3}y_{3}+F_{2}(z_{3}))x_{3}+F_{1}(z_{3})),z_{3})$\\
\item $\mathfrak{b}=(x_{3}y_{3},x_{3}y_{3}((y_{3}+F_{2}(z_{3}))x_{3}y_{3}+F_{1}(z_{3})),z_{3})$\\
\item $\mathfrak{c}=(x_{3}y_{3}(x_{3}+G_{2}(z_{3})),x_{3}y_{3}(x_{3}+G_{2}(z_{3}))(x_{3}y_{3}+F_{1}(z_{3})),z_{3})$\\
\item $\mathfrak{d}=(y_{3}(y_{3}x_{3}+G_{2}(z_{3})),y_{3}(y_{3}x_{3}+G_{2}(z_{3}))(y_{3}+F_{1}(z_{3})),z_{3}).$\\
\vspace*{5 mm}
The ideals $\mathfrak{a}$ and $\mathfrak{b}$ can be reduced in the desired form, i.e.:
\begin{equation*}
 \mathfrak{a}=(x_{3},z_{3})~\text{and}~\mathfrak{b}=(x_{3}y_{3},z_{3}).
\end{equation*}
While for the last two ideals one obtains:
\item $\mathfrak{c}=(x_{3}y_{3}(x_{3}+G_{2}(z_{3})),z_{3})=(x_{3}^{2}y_{3}+x_{3}y_{3}G_{2}(z_{3}),z_{3})$\\
\item $\mathfrak{d}=(y_{3}(y_{3}x_{3}+G_{2}(z_{3})),z_{3})=(y_{3}^{2}x_{3}+y_{3}G_{2}(z_{3}),z_{3})$.\\
\vspace*{5 mm}
Let's study $\mathfrak{c}$, in particular notice that $H_{3}(x_{3},y_{3},z_{3}):=x_{3}y_{3}\tilde{G_{2}}(z_{3})z_{3}^{i}$
is an holomorphic function which is an element of $(z_{3})$ where $\tilde{G_{2}}(0)\neq0$ (here we are using the fact that $\mu_{2}:Z_{2}\twoheadrightarrow Z_{1}$ and $G_{2}(0)=0$).
Thus
\begin{equation*}
 x_{3}^{2}y_{3}=x_{3}^{2}y_{3}+x_{3}y_{3}G_{2}(z_{3})-H_{3}(x_{3},y_{3},z_{3}),
\end{equation*}
so $\mathfrak{c}=(x_{3}^{2}y_{3},z_{3})$.\\
Repeating the same argument  for the case of $\mathfrak{d}$ one obtains $\mathfrak{d}=(y_{3}^{2}x_{3},z_{3})$.
\\
Now let's suppose the statement true for all $j=1,\dots, s-1$. Let $\tilde{\mu}=\mu_{1}\circ\cdots\circ\mu_{s}$ be a composition of
$s$ blow-ups, then $\mu_{s}$ has two possible expressions in local coordinates $(x_{s},y_{s},z_{s})$ which depend on the coordinates on $X_{s-1}$, but
by induction the local expression of $\mu_{1}\circ\cdots\circ\mu_{s-1}$ is of the form $(x_{s-1}^{i}y_{s-1}^{j},z_{s-1})$
so composing with local expressions of $\mu_{s}$ one obtains the statement.\hfill$\square$\\
\vspace*{5 mm}
Now let's recall Howald's theorem \cite{How1} for multiplier ideals of monomial ideals.
\begin{defn}
 Let $\mathfrak{a}\subset\mathbb{C}[x_{1},\dots,x_{n}]$ be a monomial ideal. We will regard $\mathfrak{a}$ as a subset of the lattice $L=\mathbb{N}^{n}$
of monomials. The \textbf{Newton Polygon} P of $\mathfrak{a}$ is the convex hull of this subset of $L$, considered as a subset of $L\otimes\mathbb{R}=\mathbb{R}^{n}$.
It is an unbounded region.
\end{defn}
\begin{nota}
 We write \textbf{1} for the vector $(1,1,\dots,1)$, which is identified with the monomial $x_{1}x_{2}\cdots x_{n}$. We use Greek letters $(\lambda\in L)$
for elements of $L$ or $L\otimes\mathbb{R}$, and exponent notation $x^{\lambda}$ for the associated monomial. For any subset $P$ of $L\otimes\mathbb{R}$,
we define $rP$ ``pointwise``
\begin{equation*}
 rP=\{r\lambda~:~\lambda\in P\}.
\end{equation*}
We write $Int(P)$ for the topoogical interior of $P$.
\end{nota}
\begin{thm}[\cite{How1}](Howald's theorem.)
 Let $\mathfrak{a}\subset\mathcal{O}_{\mathcal{A}^{n}}$ be a monomial ideal. Let $P$ be its Newton polygon. Then $\mathcal{J}(r\cdot\mathfrak{a})$ is a monomial
ideal, and contains exactly the following monomials:
\begin{equation*}
\mathcal{J}(r\cdot\mathfrak{a})=\{x^{\lambda}~:~\lambda+\textbf{1}\in Int(P)\cap L\}.
\end{equation*}
where $\mathcal{J}(r\cdot\mathfrak{a})$ is the multiplier ideal associated to $r$ and $\mathfrak{a}$ (Cfr. \cite{Laz}).
\end{thm}

Let $\mathfrak{a}$ be a monomial ideal and let $P$ be its Newton polygon. The log canonical threshold $c(\mathfrak{a})$ of $\mathfrak{a}$ is defined to be
\begin{equation}
 c(\mathfrak{a})=\sup\{r~:~\mathcal{J}(r\cdot\mathfrak{a})\neq\mathcal{O}_{X}\}.
\end{equation}
Howald's theorem shows that this must be equal to $\sup\{r~:~\textbf{1}\notin rP\}$. Thus the log canonical threshold is the reciprocal of the (unique)
number $m$ such that the boundary of $P$ intersects the diagonal in $\mathbb{R}^{n}$ at the point $m$\textbf{1}. In other words, in order to calculate
the threshold, we need only find where $P$ intersects the diagonal.
\begin{prop}[\cite{DK}]\label{pro6}
 Let $X,Y$ be complex manifolds of respective dimension $n,m$, let $\mathcal{I}\subset\mathcal{O}_{X}$, $\mathcal{J}\subset\mathcal{O}_{Y}$ be coherent ideals, and
let $K\subset X$, $L\subset Y$ be compact set. Put $\mathcal{I}\oplus\mathcal{J}:=pr_{1}^{*}\mathcal{I}+pr_{2}^{*}\mathcal{J}\subset\mathcal{O}_{X\times Y}$.
Then
\begin{equation*}
 c_{K\times L}(\mathcal{I}\oplus\mathcal{J})=c_{K}(\mathcal{I})+c_{L}(\mathcal{J}).
\end{equation*}
\end{prop}
\begin{prop}\label{pro7}
 Let $\mathcal{I}=(x^{h}y^{k},z)$ be an ideal as (\ref{pro5}) then $c(\mathcal{I})>1$.
\end{prop}
\text{Proof.} Using Howald's theorem and (\ref{pro6})
\begin{equation}
 c_{K}(\mathcal{I})=c_{K}(x^{h}y^{k})+c_{K}(z)=\dfrac{h+k}{hk}+1>1.
\end{equation}
\hfill$\square$

\subsection{Proof of Proposition (\ref{pro4})}

Now we can give the proof for the main proposition \\
\flushleft
\textit{Proof of (\ref{pro4}).} By (\ref{pro5}) and (\ref{pro6})   one can define $\displaystyle \varphi:=\text{log}(|x_{s}^{h}y_{s}^{k}|+|z_{s}|)$ thus
$c_{K}(\varphi)=c_{K}((x_{s}^{h}y_{s}^{k},z_{s}))>1$.
 Then by (\ref{pro7}) and by the definition of Lelong number the following holds:
\begin{equation*}
 \nu(\tilde{\mu_{*}}(\tilde{\Omega}),0)=
\lim_{r\rightarrow0}\dfrac{1}{r^{2}}\int_{\tilde{\mu}^{-1}(B(0,r))}\tilde{\Omega}\wedge\tilde{\mu^{*}}\left(\dfrac{i}{2\pi}\partial\bar{\partial}|z|^{2}\right)
\leq\lim_{r\rightarrow0}C_{1}r^{2(c_{K}(\varphi)-1)}= 0,
\end{equation*}
where $C_{1}$ is a positive constant.\hfill$\square$
\newpage

\subsection{Relation between mobile intersection and positive product in dimension 3}

By using Proposition (\ref{pro4}) we obtain a similar decomposition as in Theorem (\ref{pro3}), using the definition of mobile product via modifications. In fact we have the following
\begin{prop}\label{pro11}
 Let $X$ be a compact K$\ddot{\text{a}}$hler manifold of complex dimension 3 and let $\alpha_{1},\alpha_{2}\in\mathcal{M}^{2}$ then
\begin{equation*}
 \displaystyle\alpha_{1}\alpha_{2}-\langle\alpha_{1}\alpha_{2}\rangle=\left\{\sum_{t=1}^{+\infty}\nu_{\infty,t}[Y_{\infty,t}]\right\}
\end{equation*}
where $\nu_{\infty,t}\geq0$ and $Y_{t}$ are the irreducible components of codimension 2 of $\cup_{i=1}^{2}E_{nn}(\alpha_{i})$.
\end{prop}
\textit{Proof.} For $i=1,2$ Let $T_{k,i}\in\alpha_{i}$ be a sequence of currents with analytic singularities as in the Theorem (\ref{thm2}) such that $T_{k,i}\geq-\varepsilon_{k,i}\omega$ and $T_{k,i}\rightharpoonup T_{min,i}$
where $T_{min,i}$ is a positive current with analytic singularities in $\alpha_{i}.$ Now let
$\mu_{k}:X_{k}\to X$
be a common modification for $T_{k,i}$ such that $\mu_{k}^{*}T_{k,i}=[E_{k,i}]+\beta_{k,i}$, where $\beta_{k,i}$ are smooth.
Now we have that
\begin{equation*}
 T_{k,1}\wedge T_{k,2}-(\mu_{k})_{*}(\beta_{k,1}\wedge\beta_{k,2})=0~\text{on}~X\setminus (\cup_{i=1}^{2}V(\mathcal{I}(k\varphi_{min,i}))
\end{equation*}
thus by using the support theorem for currents and the inclusion $\cup_{i=1}^{2}V(\mathcal{I}(k\varphi_{min,i}))\subseteq\cup_{i=1}^{2} E_{nn}(\alpha_{i})$ we obtain
\begin{equation}\label{eq10}
 T_{k,1}\wedge T_{k,2}-(\mu_{k})_{*}(\beta_{k,1}\wedge\beta_{k,2})=\sum_{t=1}^{+\infty}\nu_{k,t}[Y_{k,t}]
\end{equation}
where $Y_{k,t}$ are the irreducible components of $E_{nn}(\alpha_{1})\cup E_{nn}(\alpha_{2})$ and $\nu_{k,t}\in\mathbb{R}$. Now let $\nu_{k,t}^{+}$ and $\nu_{k,t}^{-}$ be the positive and the negative coefficients in
the series of (\ref{eq10}), thus we have
\begin{equation}\label{eq11}
 T_{k,1}\wedge T_{k,2}+\sum_{t}(-\nu_{k,t}^{-})[Y_{k,t}]=(\mu_{k})_{*}(\beta_{k,1}\wedge\beta_{k,2})+\sum_{t}\nu_{k,t}^{+}[Y_{k,t}].
\end{equation}
 If there exists $t_{0}$ such that $\nu_{k,t_{0}}^{-}<0$ i.e. $-\nu_{k,t_{0}}^{-}>0$ then by calculating the generic Lelong number along $Y_{k,t_{0}}$, 
using  equality (\ref{eq11}), we find 
\begin{equation*}
 \nu(T_{k,1}\wedge T_{k,2},Y_{k,t_{0}})-\nu_{k,t_{0}}^{-}=\nu((\mu_{k})_{*}(\beta_{k,1}\wedge\beta_{k,2}),Y_{k,t_{0}})+0
\end{equation*}
and thanks to Proposition (\ref{pro4}) $\nu((\mu_{k})_{*}(\beta_{k,1}\wedge\beta_{k,2}),Y_{k,t_{0}})=0$ then we have
$\nu(T_{k,1}\wedge T_{k,2},Y_{k,t_{0}})=\nu_{k,t_{0}}^{-}<0$, a contradiction. Hence we have that
$\nu_{k,t}=\nu(T_{k,1}\wedge T_{k,2},Y_{k,t})$ for all $t$.
By using Lemma (\ref{lem3}) the sequence $\{\nu_{k,t}\}$ is bounded, and $T_{k,1}\wedge T_{k,2}$ are bounded in mass since they are in
the same cohomology class for all $k$, while $\{(\mu_{k})_{*}(\beta_{k,1}\wedge\beta_{k,2})\}$ are bounded in mass thanks to Theorem (\ref{thm3}). 
Hence we can extract  a common subsequence such that  passing
in cohomology we obtain the statement.\hfill$\square$ 

\newpage
And as a consequence we have

\begin{cor}
 Assume $\dim(X)=3$ and $\alpha,\beta\in\mathring{\mathcal{M}}^{2}$. Let $T_{min,\alpha},T_{min,\beta}$ be positive closed currents with
minimal singularities in $\alpha$ and $\beta$ respectively. Then for all $x\in X$ one has
\begin{equation}
 \lim_{k\to+\infty}\nu(T_{k,\alpha}\wedge T_{k,\beta},x)=\nu(T_{min,\alpha}\wedge T_{min,\beta},x),
\end{equation}
where $\{T_{k,\alpha}\}$ and $\{T_{k,\beta}\}$ are sequences of currents as in Theorem (\ref{thm2}) which weakly converge
to $T_{min,\alpha}$ and $T_{min,\beta}$ respectively.
\end{cor}

\textit{Proof.} By using Theorem (\ref{pro3}) we get the following decomposition, i.e. the Siu decomposition of $T_{min,\alpha}\wedge T_{min,\beta}$:
\begin{equation}
 T_{min,\alpha}\wedge T_{min,\beta}-\langle T_{min,\alpha}\wedge T_{min,\beta}\rangle=\sum_{h=1}^{N}\nu(T_{min,\alpha}\wedge T_{min,\beta},Y_{h})[Y_h]
\end{equation}

On the other hand by Proposition (\ref{pro11}) we also have the following decomposition
\begin{equation}
 T_{k,1}\wedge T_{k,2}-(\mu_{k})_{*}(\beta_{k,1}\wedge\beta_{k,2})=\sum_{h=1}^{N}\nu(T_{k,1}\wedge T_{k,2},Y_{h})[Y_{h}]
\end{equation}
where in both decompositions $Y_{h}$ are the irreducible components of dimension $1$ of $E_{nn}(\alpha)\cap E_{nn}(\beta)$.\\
Then by choosing a common subsequence we obtain by uniform bound of mass and by Proposition (\ref{pro9})  that 
$\nu(T_{k,1}\wedge T_{k,2},Y_{h,k})\to\nu_{\infty,h}\leq\nu(T_{min,\alpha}\wedge T_{min,\beta},Y_{h})$.
Hence we have  positive current
\begin{equation}
\displaystyle S=T_{min,\alpha}\wedge T_{min,\beta}-\langle T_{min,\alpha}\wedge T_{min,\beta}\rangle-\sum_{h=1}^{N}\nu_{\infty,h}[Y_{h}]
\end{equation}
 whose cohomology class is the class zero therefore $S=0$ then the statement holds.\\
We can also conclude that
\begin{equation}
 (\mu_{k})_{*}(\beta_{k,1}\wedge\beta_{k,2})\rightharpoonup\langle T_{min,\alpha}\wedge T_{min,\beta}\rangle.
\end{equation}

\hfill$\square$

\cleardoublepage
\bibliography{Bibliography}
\bibliographystyle{amsalpha}
\end{document}